\documentclass[11pt]{amsart}
\usepackage{amsmath,amsfonts,amssymb,amsthm,amscd,hyperref}
\usepackage[english]{babel}
\usepackage[applemac]{inputenc}
\usepackage[mathscr]{eucal}
\usepackage{graphicx,url}
\usepackage{appendix}
\usepackage{color}
\newtheorem{theorem}{Theorem}[section]
\newtheorem{prop}[theorem]{Proposition}

\newtheorem{remark}[theorem]{Remark}
\theoremstyle{definition}

\newcommand{\R}{\mathbb{R}}

\title[]{Collisions of vortex filament pairs}
\author[V. Banica]{Valeria Banica}
\address[V. Banica]
{Laboratoire de Math\'ematiques et de Mod\'elisation d'\'Evry (UMR 8071)\\ 
Universit\'e d'\'Evry, 23 Bd. de France, 91037 Evry\\ 
France, Valeria.Banica@univ-evry.fr} 

\author[E. Faou]{Erwan Faou}
\address[E. Faou]
{INRIA \& ENS Cachan Bretagne  \\
Avenue Robert Schumann F-35170 Bruz, 
France, Erwan.Faou@inria.fr}

\author[E. Miot]{Evelyne Miot}
\address[E. Miot]
{Centre de Math\'ematiques Laurent Schwartz (UMR 7640) \\ 
\'Ecole Polytechnique 91128 Palaiseau\\ 
France, Evelyne.Miot@math.polytechnique.fr} 

\thanks{First and last authors are partially supported by the French ANR project SchEq ANR-12-JS-0005-01. The second author is supported by the ERC starting grant GEOPARDI No. 279389. The last author is partially supported by the French ANR project GEODISP ANR-12-BS01-0015-01.}

\begin{document}
\maketitle 

\begin{abstract}
We consider the problem of collisions of vortex filaments  for a model  introduced by Klein, Majda and Damodaran  \cite{KlMaDa} and Zakharov \cite{Zh1,Zh2} to describe the interaction of almost parallel vortex filaments in three-dimensional fluids. Since the results of Crow \cite{Cr} examples of collisions are searched as perturbations of antiparallel translating pairs of filaments, with initial perturbations related to the unstable mode of the linearized problem; most results are numerical calculations. In this article we first consider a related model for the evolution of pairs of filaments and we display another type of initial perturbation leading to collision in finite time.  Moreover we give numerical evidence that it also leads to collision through the initial model. We finally study the self-similar solutions of the model. 
\end{abstract}
\tableofcontents

\section{Introduction}
In this article we investigate the possible collision scenarios for a system of equations modelling the dynamics of vortex filaments in three-dimensional incompressible fluids. Since the results of Crow \cite{Cr},  collisions in this model  are searched as perturbations of an antiparallel pair of vortex filaments, with initial perturbations related to the unstable mode of the linearized problem; most results are numerical calculations.

Let us now present the model at stake and survey the known results on it. 
 {It was} introduced by Klein, Majda and Damodaran  \cite{KlMaDa} and Zakharov \cite{Zh1,Zh2} in order to describe the evolution of $N$ almost parallel vortex filaments in a three-dimensional fluid governed by the Navier-Stokes equations.  According to this model, the vortex filaments are curves parametrized by 
 $$(\Re \psi_j(t,\sigma), \Im \psi_j(t,\sigma),\sigma),\quad \sigma\in \R,\quad 1\leq j\leq N,$$
where $\Psi_j:\R\times \R\to \mathbb{C}$, and their dynamics is given by the following Hamiltonian system of one-dimensional Schr\"odinger  equations with vortex type interaction
\begin{equation}\label{KMD}
\begin{cases}
\displaystyle 
i\partial_t \Psi_j+\alpha_j\Gamma_j \partial_\sigma^2\Psi_j+\sum_{k\neq j} \Gamma_k \frac{\Psi_j-\Psi_k}{|\Psi_j-\Psi_k|^2}=0,\quad 1\leq j\leq N,\\
\Psi_j(0,\sigma)=\Psi_{j,0}(\sigma). 
\end{cases}
\end{equation}
Here $\alpha_j\in\mathbb R$ is a parameter related to the core structure  and $\Gamma_j$ is the circulation of the $j$th filament. Particular solutions are the parallel filaments $(X_j)$ with $X_j(t,\sigma)=(z_j(t),\sigma)$ for $(z_j)$ solutions of the 2-D point vortex system
\begin{equation}\label{point-vortex}
\displaystyle 
i\dot{z}_j+\sum_{k\neq j} \Gamma_k \frac{z_j-z_k}{|z_j-z_k|^2}=0,\quad 1\leq j\leq N.
\end{equation}
In this paper we will focus on pairs of filaments, namely $N=2$.

\medskip

The system \eqref{KMD} appears as the leading-order asymptotic equation under particular conditions relating the wavelength, the amplitude of the perturbations around $(X_j)$, the core thicknesses of the filaments, and the distances between the filaments. The attention was focused {in \cite{KlMaDa,Zh1,Zh2}} on perturbations of two parallel filaments $(X_j)$ with the same core parameters, $\alpha_1=\alpha_2$.  In the case of circulations having the same sign, 
the perturbations around the corotating pair of parallel filaments $(X_j)$ were proved to be linearly stable and global existence was conjectured, while for 
circulations of opposite signs instability was shown to appear for the linearized operator and numerical evidence of finite-time collapse was given.  We send the reader to chapter 7 in \cite{MaBe} for a complete description of the model \eqref{KMD} and of these results.

Later Kenig, Ponce and Vega  \cite{KePoVe} proved global existence for \eqref{KMD} with same core parameters and same circulations for perturbations around a corotating pair of parallel filaments $(X_j)$ and for perturbations around $(X_j)$ with $(z_j)$ a rotating equilateral triangle.

Recently in \cite{BaMi, BaMi2} symmetric perturbations around  $(X_j)$, with $(z_j)$ a rotating $N-$polygon,  were proved to exist globally in time for \eqref{KMD} with same core parameters and circulations. Moreover, existence of travelling waves was displayed for \eqref{KMD}, together with collision scenarios for at least three filaments.

Finally, let us mention  two other related results. Lions and Majda \cite{LiMa} established global existence without uniqueness of weak solutions of \eqref{KMD} as perturbations around any $(X_j)$ with same core parameters and circulations. A work by Craig and Garcia \cite{CrGa} is in progress on quasi-periodic and periodic perturbations of corotating pairs of parallel filaments using KAM theory. 

Summarizing, for pairs of filaments no collision for \eqref{KMD} has been displayed rigorously. Moreover, the only numerical simulations of collisions have been performed starting from a perturbation given by the unstable mode of the linearized system around the antiparallel translating pair of filaments (given by the solution of the 2D point vortex system \eqref{point-vortex} with $\Gamma_1=-\Gamma_2$). In this paper we will  construct exact examples of collisions of pairs of filaments in related models together with numerical simulations. More precisely, in Section \S \ref{sect-opposite} we consider the system obtained from \eqref{KMD} by setting $\alpha_1=-\alpha_2$ and we 
 display another type of explicit initial perturbation {of the antiparallel translating pair of filaments} leading to an exact collision for this system. Moreover we give numerical evidence that this initial perturbation also leads to collision in \eqref{KMD} with $\alpha_1=\alpha_2$. Remarkably, this type of data is precisely the one used in \cite{BaMi} to construct an explicit collision in \eqref{KMD} for any number of filaments blue larger or equal to three. Then in Section {\S \ref{sect-selfsim}} we establish the existence of self-similar solutions of \eqref{KMD} with $\alpha_1=\alpha_2$. Finally, the last Section \S \ref{secnum} gathers all the numerical investigations.

\section{Antiparallel vortex filament pairs collisions}\label{sect-opposite}
We consider opposite circulations
$$\Gamma_1=-\Gamma_2=1.$$
Then, for two nearly parallel vortex filaments, system \eqref{KMD} writes
\begin{equation}\label{KMD2gen}
\begin{cases}
\displaystyle i\partial_t \Psi_1+\alpha_1\partial_\sigma^2\Psi_1-\frac{\Psi_1-\Psi_2}{|\Psi_1-\Psi_2|^2}=0,\\
\displaystyle i\partial_t \Psi_2-\alpha_2\partial_\sigma^2\Psi_2-\frac{\Psi_1-\Psi_2}{|\Psi_1-\Psi_2|^2}=0.\\
\end{cases}
\end{equation}

Moreover,
$$z_1(t)=-\frac {it}2+1,\quad z_2(t)=-\frac{it}2-1,$$
is an explicit solution of the vortex point system
\begin{equation*}
\begin{cases}
\displaystyle i\partial_t z_1- \frac{z_1-z_2}{|z_1-z_2|^2}=0,\\
\displaystyle i\partial_t z_2- \frac{z_1-z_2}{|z_1-z_2|^2}=0.
\end{cases}
\end{equation*}

\subsection{Case of opposite core parameters}

We start by studying the system \eqref{KMD2gen} with $$\alpha_1=-\alpha_2=1.$$ We consider
the Gaussian 
$$G(\sigma)=\frac{e^{-\frac{\sigma^2}{1-4i}}}{\sqrt{1-4i}}.$$ 
Below we construct a global solution of \eqref{KMD2gen} such that the filaments $\psi_1(t,\sigma)$ and $\psi_2(t,\sigma)$ collapse at time $t=1$ and height $\sigma=0$. 
\begin{prop}
\label{prop21}
The perturbation $\Psi_{j}(t,\sigma)$ of $X_j(t,\sigma)=(z_j(t),\sigma)$ generated by the data
$$\Psi_{1,0}(\sigma)=1-G(\sigma),\quad \Psi_{2,0}(\sigma)=-1+G(\sigma)$$
evolves through a global solution of system \eqref{KMD2gen} with $\alpha_1=-\alpha_2=1$ in the following (mild) sense:
\begin{equation*}
\begin{split}
\Psi_1(t,\sigma)&=-\frac{it}2+e^{it\partial_\sigma^2}(1-G)(\sigma)-\frac 12D(t,\sigma),\quad t\in \R,\\
 \Psi_2(t,\sigma)&=-\frac{it}2-e^{it\partial_\sigma^2}(1-G)(\sigma)-\frac 12D(t,\sigma),\quad t\in \R,
\end{split}
\end{equation*}
where $$ D\in C (\R, L^q(\R)), \quad \forall q\in[2,+\infty)$$
and is given by
$$D(t,\sigma)=i\int_0^t e^{i(t-s)\partial_\sigma^2}\left(\frac{1}{1-\overline{e^{is\partial_\sigma^2}G(\sigma)}}-1\right)\,ds,\quad t\in \R.$$

In particular the filaments $\psi_1(t,\sigma)$ and $\psi_2(t,\sigma)$ collide at time $t=1$ and height $\sigma=0$. 

\end{prop}
\begin{remark}
We notice that the $\Psi_j(t,\sigma)$ are symmetric perturbations of $X_j(t,\sigma)$ plus a non-linear shift $D$ that acts for all times $t$ and heights $\sigma$, so the configuration does not keep any symmetry (see Figure \ref{opposite} in Section \S\ref{secnum}). 

\end{remark}
\begin{remark}
Once we have an initial configuration $\psi_{j,0}$ leading to collision, any shifted perturbation $\psi_{j,0}(\sigma)+u_0(\sigma)$, with $u_0\in L^2(\mathbb R)$, gives rise to the solution $\Psi_j(t,\sigma)+e^{it\partial_\sigma^2}u_0(\sigma)$ and therefore  also leads to collision at the same time and height.
\end{remark}

\begin{proof}
Adding and substracting the equations in \eqref{KMD2gen} we get  a system of equations for the functions $\Psi=\Psi_1-\Psi_2$ and $\Phi=\Psi_1+\Psi_2:$
\begin{equation}\label{KMD2diff}
\begin{cases}
\displaystyle i\partial_t \Psi+\partial_\sigma^2 \Psi=0\\
\displaystyle i\partial_t \Phi+\partial_\sigma^2\Phi=2\frac{\Psi}{|\Psi|^2}.
\end{cases}
\end{equation}
The solution of the first equation is
\begin{equation*}
\begin{split}
\Psi(t,\sigma)&=e^{it\partial_\sigma^2}\Psi_{0}(\sigma),\quad (t,\sigma)\in \R\times \R.\end{split}\end{equation*}
We look then for a mild solution of the second equation, namely
\begin{equation*}
\begin{split}
\Phi(t,\sigma)&=e^{it\partial_\sigma^2}\Phi_0(\sigma)-2i\int_0^t e^{i(t-s)\partial_\sigma^2}\frac{\Psi}{|\Psi|^2}(s,\sigma)\,ds.\end{split}
\end{equation*}
Since we have 
$$\Psi_0(\sigma)=2(1-G(\sigma)),$$
we retrieve the explicit example of  \cite{BaMi}:\begin{equation*}
\begin{split}\Psi(t,\sigma)&=2e^{it\partial_\sigma^2}(1-G)(\sigma)=2\left(1-\frac{e^{-\frac{\sigma^2}{1-4i(1-t)}}}{\sqrt{1-4i(1-t)}}\right),\quad (t,\sigma)\in \R\times \R,
\end{split}\end{equation*}
so that collision occurs if and only if $(t,\sigma)=(1,0)$, 
$$\Psi(1,0)=e^{i1\partial_\sigma^2}\Psi_0(0)=0.$$
Note that for simplicity we have also fixed the initial symmetry
$$\Phi_0(\sigma)=0,$$
so 
$$\Phi(t,\sigma)=-i\int_0^t e^{i(t-s)\partial_\sigma^2}\frac{1}{1-\overline{e^{is\partial_\sigma^2}G(\sigma)}}\,ds.$$
We obtain then the explicit expressions of $\Psi_1(t,\sigma)$ and $\Psi_2(t,\sigma)$ in the statement of the Proposition.

Finally we justify that $D\in C(\R,L^q(\R))$, thus $\Phi$ is globally defined even beyond the collision time $t=1$. First we notice that there exists $m>0,C>0$ such that
\begin{equation}\label{estG}\left|\frac{1}{1-e^{is\partial_\sigma^2}G(\sigma)}\right|\leq C\max\left\{1,\frac 1{\sigma^2+|1-s|}\right\},\quad \forall |\sigma|\leq m,\quad  s\in \R.
\end{equation}
Indeed, by a second order Taylor expansion in variables $(s,\sigma)$ for the function $$f(s,\sigma)=e^{is\partial_\sigma^2}G(\sigma)=\frac{e^{-\frac{\sigma^2}{1-4i(1-s)}}}{\sqrt{1-4i(1-s)}}$$ we obtain the existence of $m>0$ and $0<\delta_0<1$ such that
\begin{equation*}
\begin{split}
f(1,0)-f(s,\sigma)&=-2i(1-s)-\sigma^2-6(1-s)^2+r(s,\sigma)
,\end{split}\end{equation*}
with
$$|r(s,\sigma)|\leq \frac {\sigma^2+(1-s)^2}{2}\leq \frac {\sigma^2+|1-s|}{2},\quad \forall |\sigma|\leq m, \quad  |s-1|\leq \delta_0.$$
So if $|\sigma|\leq m$ and $|s-1|\leq  \delta_0$,
$$\frac 1{|f(1,0)-f(s,\sigma)|}\leq \frac 1{|\sigma^2+6(1-s)^2+2i(1-s)|-|r(s,\sigma)|}\leq \frac2{\sigma^2+|1-s|}.$$
On the other hand for $|s-1|\geq \delta_0$ we have $|f(s,\sigma)|\leq (1+16\delta_0^2)^{-1/4}$, so we get the existence of $C>0$ such that 
 $$\frac 1{|f(1,0)-f(s,\sigma)|}\leq C.$$
This establishes \eqref{estG}.  We infer that there exists $C>0$ such that
\begin{equation}\label{ineq:f}\left|\frac{e^{is\partial_\sigma^2}G(\sigma)}{1-e^{is\partial_\sigma^2}G(\sigma)}\right|\leq  C\frac{e^{-\frac{\sigma^2}{1+16(1-s)^2}}}{{ (1+16(1-s)^2)^\frac 14}}\max\left\{1,\frac 1{\sigma^2+|1-s|}\right\},\quad \forall  (\sigma,s)\in\R\times \R.
\end{equation}

Now we consider $2\leq q<\infty$. For $t>0$, using {dispersion estimate we get
\begin{equation*}
\|D (t)\|_{L^q}\leq C
 \int_{0}^{t}\left\|\frac{e^{is\partial_\sigma^2}G(\sigma)}{1-e^{is\partial_\sigma^2}G(\sigma)}\right\|_{L^{\overline{q}}}\frac {ds}{|t-s|^{\frac 12-\frac 1q}},
 \end{equation*}
  where $\overline{q}$ denotes the conjugate exponent of $q$.
We split the integral in $\sigma$ into the regions $|\sigma|\leq
1,1\leq|\sigma|\leq\max\{1,|1-s|^2\},\max\{1,|1-s|^2\}\leq|\sigma|$ and use
\eqref{ineq:f} to get
\begin{equation*}
\left\|\frac{e^{is\partial_\sigma^2}G(\sigma)}{1-e^{is\partial_\sigma^2}G(\sigma)}\right\|_{L^{\overline{q}}}\leq
\left\|1+\frac{1}{\sigma^2+|1-s|}\right\|_{L^{\overline{q}}(|\sigma|\leq
1)}+\left\|\frac{e^{-\frac{\sigma^2}{17|1-s|^2}}}{|1-s|^\frac
12}\right\|_{L^{\overline{q}}}+\left\|e^{-\frac{\sigma}{17}}\right\|_{L^{\overline{q}}}.
\end{equation*}
We perform in the integrals above changes of variables and we obtain
\begin{equation*}
\|D (t)\|_{L^q}\leq
C\int_0^t\left(1+|1-s|^{-\frac12-\frac 1{2q}}+|1-s|^{\frac 12-\frac
1q}\right)\frac {ds}{|t-s|^{\frac 12-\frac 1q}}\leq C(1+t).
\end{equation*}}
When $t<0$ the argument goes the same. 
This proves that $D\in L^\infty_{\text{loc}}(\R, L^q(\R))$. The time continuity can be established by the same arguments.
 \end{proof}

\subsection{Case of positive core parameters}

In the previous subsection we have obtained an initial configuration that leads to collisions in the case when the numbers $\alpha_1,\alpha_2$ have opposite signs. However, in the derivation of \eqref{KMD} these numbers are vortex core parameters supposed to be positive. 
We consider again opposite circulations but same sign core parameters 
$$\Gamma_1=-\Gamma_2=1,\,\,\,\alpha_1=\alpha_2=1.$$
So system \eqref{KMD2gen} is
\begin{equation}\label{KMD2same}
\begin{cases}
\displaystyle i\partial_t \Psi_1+\partial_\sigma^2\Psi_1-\frac{\Psi_1-\Psi_2}{|\Psi_1-\Psi_2|^2}=0,\\
\displaystyle i\partial_t \Psi_2-\partial_\sigma^2\Psi_2-\frac{\Psi_1-\Psi_2}{|\Psi_1-\Psi_2|^2}=0.\\
\end{cases}
\end{equation}
 As we will see in the Section devoted to numerical simulations, it is worth mentioning that the same initial data from the previous subsection $\Psi_1(0)=1-G$,  $\Psi_2(0)=-1+G$ yields a very nice collision behavior, see Figure \ref{samenosym}  in Section \S\ref{secnum}.

\begin{remark}
To establish the occurence of collision we need a solution of \eqref{KMD2same} such that $(\Psi_1-\Psi_2)(t,\sigma)$ values $2(1-G(\sigma))$ at time $t=0$, stays close to $2(1-e^{it\partial_\sigma^2}G)$ and vanishes at $(t,\sigma)=(1,0)$ . 
The difference $\Psi_1-\Psi_2$ solves the equation
\begin{equation}\label{general}
\partial_t^2f+\partial_\sigma^4 f+\frac{2\partial_\sigma^2 \overline{f}}{\overline{f}^2}-\frac{4(\partial_\sigma\overline{f})^2}{\overline{f}^3}=0.
\end{equation}
It is not obvious to perform a perturbative argument for constructing such solutions - note for instance that the first nonlinear term valued on the ansatz $2(1-e^{it\partial_\sigma^2}G)$ is not in $L^1([0,1],L^2)$. 
\end{remark}
\begin{remark}One can write the difference $(\Psi_1-\Psi_2)(t,\sigma)$ using an integral equation,
\begin{equation}\label{eqintdiff}\begin{split}
&\int\left(\frac{e^{i\frac{(\sigma-y)^2}{4t}}}{\sqrt{4\pi it}}\Psi_{1,0}(y)-\frac{e^{-i\frac{(\sigma+y)^2}{4t}}}{\sqrt{4\pi it}}\Psi_{2,0}(y)\right)dy\\
&+2\int_0^t\int\frac{\sin\frac{\sigma^2+y^2}{4(t-\tau)}e^{-i\frac{\sigma y}{2(t-\tau)}}}{\sqrt{4\pi i(t-\tau)}}\frac 1{\overline{\Psi_1-\Psi_2}(\tau,y)}dy d\tau,\end{split}
 \end{equation}
 but it is not obvious to control pointwise the second term in order to get cancellation for some $t>0$ at some point $\sigma$. 
\end{remark}

Next we consider as in \cite{Zh1,Zh2,MaBe} an initial configuration of two filaments satisfying
$$\Psi_1(0)=-\overline{\Psi}_2(0).$$
Then, as long as the solution $(\Psi_1,\Psi_2)$ to \eqref{KMD2gen} exists, it satisfies
$$\Psi_1(t)=-\overline{\Psi}_2(t),$$
because $(-\overline{\Psi}_2,-\overline{\Psi}_1)$ is also solution to \eqref{KMD} with same initial datum.
In particular $\Psi_1-\Psi_2=2\Re(\Psi_1)$ is real and the equation for $\Psi_1$ in \eqref{KMD} reduces to
\begin{equation}\label{eqpsi}
 i\partial_t \Psi_1+\partial_\sigma^2\Psi_1-\frac{1}{2\Re(\Psi_1)}=0.
\end{equation}
Collision occurs when the real part of $\psi_1$ vanishes. With the same initial data $\psi_1(0)=1-G$ as before, numerical simulations provide a symmetric\footnote{In the sense that $\Psi_1(t)=-\overline{\Psi}_2(t)$.} collision, see Figure \ref{same} in Section \S\ref{secnum}.

\begin{remark}Note that the source term$$
\int_0^t e^{i(t-s)\partial_\sigma^2}\frac{1}{2\Re(\Psi_1)(s)}\,ds$$
valued on the ansatz $\Psi_1(s,\sigma)= 2(1-e^{it\partial_\sigma^2}G)$ does not belong to $L^1([0,1]L^2)$, since the behavior of $\Re(1-e^{it\partial_\sigma^2}G)$ near $(t,\sigma)=(1,0)$ is $\sigma^2+6(1-t)^2$. This is in contrast with the non-symmetric data case for which the term to be  evaluated is $|1-e^{it\partial_\sigma^2}G|^{-1}\leq C(\sigma^2+|1-t|)^{-1}$ and has been handled in Proposition \ref{prop21}. This is an obstacle in performing a perturbative argument to construct solutions of \eqref{eqpsi} around $1-e^{it\partial_\sigma^2}G$.  
\end{remark}

We conclude this Section with a discussion on another approach to exhibit collision scenarios. Namely, one could try the method used by Merle and Zaag \cite{MeZa} for the problem of vortex reconnection
with the boundary in a superconductor under the planar approximation: for $\Psi_1$ solution of \eqref{eqpsi} we denote
$u(t,\sigma)=\frac 1{\Psi_1(t,\sigma)+\frac{it}{2}}$ that solves
\begin{equation}\label{GPD}i\partial_t u+ \partial_\sigma^2u-\frac{u^2}{2}\left(1-\frac{|u|^2}{\Re (u)}\right)-2\frac{(\partial_\sigma u)^2}{u}=0,\end{equation}
with boundary condition $1$ at infinity. Now we note that 
$$\Re(\Psi_1)=\frac{\Re(u)}{|u|^2},$$
so in order to have a cancellation for $\Re(\Psi_1)$ it is enough to obtain a solution $u(t,\sigma)$ blowing-up in $L^\infty$ in finite time at one point. In \cite{MeZa} instead of \eqref{GPD} a heat equation is obtained with power nonlinearity and a square gradient term. A pointwise blowing-up solution is constructed starting from the explicit profile of pointwise blowing-up solution of the heat equation with power nonlinearity that was described in \cite{MeZa2}. In here \eqref{GPD} is a nonlinear Schr\"odinger type equation and we are looking for what is called a dispersive blow-up solution. This kind of problem has been considered recently by Bona, Ponce, Saut and Sparber \cite{BPSS}. Introducing $u=v+1$, the equation on $v$ is of Gross-Pitaevskii type
\begin{equation}\label{GPDbis}i\partial_t v+ \partial_\sigma^2v+\frac 12\Re(v)+F(v)=0,\end{equation}
with
$$F(v)=\Re v\left(v+\frac{(\Re v)^2}{2}+i\Re v\Im v\right)+\frac{(\Im v)^2}{2}+i(\Im v)^3-\frac{(\Im v)^4}{1+\Re(v)}-2\frac{(\partial_\sigma v)^2}{1+v},$$
and \cite{BPSS} provides an explicit initial data that leads to dispersive blow-up for the linearized Gross-Pitaevskii equation:
$$v_0^m(\sigma)=\frac{e^{ix^2}}{(1+x^2)^m}+B_0\left(\frac{\cos(x^2)}{(1+x^2)^m}\right),\quad \frac 14<m \leq \frac 12,$$
where $B_0$ has the symbol $-\frac{1}{(1+\xi^2)(\sqrt{\xi^2(1+\xi^2)^{-1}+1)}}$.
Numerics again exhibits a collision in finite time by taking $\Psi_1(0,\sigma)=\frac 1{1+v^m_0(\sigma)}$.\\
The article \cite{BPSS} gives also a recipe to conclude the dispersive blow-up for nonlinear equations of type \eqref{GPDbis}, from the same initial data. It consists in controlling the Duhamel term, once the existence of a local nonlinear solution in some appropriate spaces is known. So up to multiplying the initial data with a small constant, a local Cauchy theory for small data regular as $v_0^m$ for \eqref{GPDbis} is sufficient for our purposes. The regularity of the initial data is proved in \cite{BPSS} to be 
$$v_0^m\in\mathcal C^\infty\cap L^\infty\cap H^s(\mathbb R),\forall s<2m-\frac 12.$$
In our case, for proving local in time wellposedness of \eqref{GPDbis}, the first terms in $F(v)$ can be treated at the $L^2$ level in $ L^\infty([0,T[,L^2)\cap L^4([0,T[,L^{\infty})$.  Then for treating $\frac{(\Im v)^4}{1+\Re v}$ one needs to control through the fixed point argument a positive lower  bound for $1+\Re v$, which is not obvious. Finally,  inglobing in the fixed point the term involving $(\partial_\sigma v)^2$ is not obvious neither, since the initial data is of regularity at most $\mathcal C^\infty\cap L^\infty\cap H^{\frac 12^-}(\mathbb R)$.


\section{Self-similar vortex filament pairs collisions}\label{sect-selfsim}

In this section we investigate the existence of a particular class of solutions to \eqref{KMD} in the case
$$\Gamma_1=-\Gamma_2=1,\quad \alpha_1=\alpha_2=1.$$ More precisely, 
as suggested in \cite{Zh1,Zh2,KlMaDa}  we look for a self-similar solution of \eqref{eqpsi}:
$$\psi_1(t,\sigma)=\sqrt{t}\:u\left(\frac{\sigma}{\sqrt{t}}\right).$$ 
Setting $x=\sigma/\sqrt{t}$, note that if the profile $u(x)$ behaves like $\alpha |x|$ when $|x|\to +\infty$  then the corresponding two filaments are localized at $+\infty$ near the lines passing through the origin and generated by $(1,\Re\alpha,\Im\alpha)$ and $(1,-\Re\alpha,\Im\alpha)$ respectively, { and at $-\infty$ near the lines passing through the origin and generated by $(1,-\Re\alpha,-\Im\alpha)$ and $(1,\Re\alpha,-\Im\alpha)$ respectively}. The filaments are not parallel, so a priori this kind of configuration does not enter the setting of the modelisation  \eqref{KMD} proposed by  \cite{KlMaDa}. As \cite{Zh1,Zh2} mention, this is actually a case for getting information on possible reconnection scenarios. 

We obtain as an equation for the profile $u(x)$ 
\begin{equation}\label{equ}
i(u-xu')+2u''-\frac{1}{\Re(u)}=0.
\end{equation}


Up to a change of scale we may rewrite Equation \eqref{equ} as
\begin{equation}\label{syst:self-sim}
\begin{cases}
\displaystyle v' = i x v - \frac{x}{ \Re(u) }\\
\displaystyle u - x u' = v.
\end{cases}
\end{equation}
We look for a solution such that
$$v \text{ is even},\quad  v(0)=1.$$ The condition $v(0)=1$ corresponds to starting  {at time $t=1$ with two filaments at distance $2$ at the origin level $\sigma=0$}.
For the function $u$ we impose the following conditions at infinity:
\begin{equation}\label{infcond}
u(x)\sim \alpha |x|,\quad |x|\to+\infty
\end{equation}
where $\alpha\in \mathbb{C}$ is such that $\Re(\alpha)>0$. 

\begin{theorem}\label{thm:ak}Let $\alpha \in \mathbb{C}$.  Let
$$E=\left\{w\in C^1(\R),\:w(0)=w'(0)=0,\:\: w\text{  even  },\:\: \|w\|_{L^\infty} + \||x|^{-1}w'\|_{L^\infty}\leq \frac{\alpha}{4}\right\}.$$
There exists a numerical constant $K_0>1$ such that if $\Re(\alpha)>K_0$, there exists a unique $v \in 1+ E$ such that the couple $(u,v)$, with $u$ defined by
\begin{equation}\label{couplage}
u(x) = 1+|x|\left( \alpha + \int_{|x|}^{+\infty} \frac{v(z)-1}{z^2} \,dz\right), \forall x\neq 0,\quad u(0)=1,
\end{equation}
is a solution of the system \eqref{syst:self-sim} satisfying the condition \eqref{infcond}. Moreover,
\begin{equation}\label{poscond}
\Re (u(x))>0, \,\forall x\in\mathbb R,
\end{equation}
and $u$ is a Lipschitz function on $\R$ with $u'\in C(\R^\ast)$.
\end{theorem}

\begin{remark}
We could replace $v-1$ by $v-v(0)$ with $v(0)$ any complex number satisfying $\Re(v(0))>0$ (so that the filaments are separated at $\sigma=0$ and $t=1$).
\end{remark}
\begin{remark}
We shall see below that if $v\in 1+E$ then $\int_0^{+\infty}|v-1|/z^2\leq\alpha /2$, so that $u'$ is continuous on $\R^\ast$ and has a jump discontinuity at $x=0$ with $u'(0)^+=-u'(0)^-=\alpha+\int_0^{+\infty}(v-1)/z^2\neq 0$. This means that the filaments exhibit a corner at $\sigma=0$.
\end{remark}

\begin{proof}
Note that since the system \eqref{syst:self-sim} is invariant by the transformation $(u,v)\mapsto (u-i\Im(\alpha)|x|,v)$ we may assume that $$\alpha\in \R, \quad \alpha>0.$$

We will show by a fixed point argument that there exists a unique solution as in Theorem \ref{thm:ak} such that $v-1$ belongs to the Banach space $E$ endowed with the metric
$$d(w_1,w_2)=\|w_1-w_2\|_{\infty}+\||x|^{-1}(w_1-w_2)\|_{\infty}.$$

\medskip

Before proving this, let us derive first its consequences.
Let $(u,v)$ be such a solution of \eqref{syst:self-sim} with $v-1\in E$. 
Integrating the second equation of \eqref{syst:self-sim} on the intervals $[x,+\infty)$ when $x>0$ or $(-\infty,x]$ if $x<0$ we find
$$
u(x)= x \alpha  + x \int_x^{+\infty} \frac{v(z)}{z^2} dz\quad \text{for  } x>0
$$
and
$$
u(x)= -x \alpha  - x \int_{-\infty}^x \frac{v(z)}{z^2} dz\quad \text{for }x<0.
$$
{In particular the first part of \eqref{couplage} is satisfied.}
	Since $v$ is even we infer that $u$ is even as well. Note that $u'$ has a jump discontinuity at $x=0$ but $u-xu'{=v}\in C^1(\R)$. From now on we only consider $x>0$.

We now check that $\Re(u)>0$ on $\R$. 
We set
$$w=v-v(0)=v-1.$$
In order to consider converging integrals we rewrite $u$ as
\begin{equation}\label{def:u}
u(x) = 1 + x \alpha  + x \int_x^{+\infty} \frac{w(z)}{z^2} \,dz
    =1+x\left( \alpha + \int_x^{+\infty} \frac{w(z)}{z^2} \,dz\right).
\end{equation}
Note that the integral is well-defined since $w$ is bounded.
We first observe that {$w\in E$ implies}
\begin{equation} \label{ineq:sup}\sup_{x\in \R_+}\left|\int_x^{+\infty} \frac{w}{z^2}\,dz\right|< \frac{\alpha}{2}.\end{equation}
Indeed, {for $w\in E$,}
\begin{equation} \label{ineq:int}
\begin{split}
\left|\int_x^{+\infty} \frac{w}{z^2}\,dz\right|&\leq 
\int_0^{1} \frac{|w|}{z^2}\,dz+\int_{1}^{+\infty} \frac{|w|}{z^2}\,dz\leq \frac{\||x|^{-1}w'\|_{\infty}}{2}+ {\|w\|_{\infty}},
\end{split}
\end{equation}
where we have used that { $w(0)=0$ and}
$$|w(z)|\leq \int_0^z |w'(t)|\,dt\leq \frac{z^2\||x|^{-1}w'\|_{\infty}}{2},$$ so the inequality \eqref{ineq:sup} follows.

In particular \eqref{def:u} and \eqref{ineq:sup} yield
\begin{equation}\label{lower-bound:u}\Re(u(x))\geq 1+\frac{x\alpha}{2}>0,\quad \forall x\geq 0,\end{equation}
{together with $u(0)=1$, so all the claims of the Theorem are verified. 

\medskip

We are left with the proof of the existence of solutions in $E$. In view of the first equation in \eqref{syst:self-sim}} we define the following operator:
\begin{equation*}\begin{split}
P(w)(x)= e^{\frac{ix^2}{2}}-1-e^{\frac{ix^2}{2}}\int_0^x \frac{y e^{-\frac{iy^2}{2}}}{\Re(u(y))}\,dy,\quad x\in \R,\end{split}\end{equation*} 
{where $u(y)$ is defined by \eqref{def:u}.}

We next show that $P$ has a unique fixed point $w\in E$ if $\alpha>K_0$ with $K_0$ sufficiently large.

When $w\in E$ the function $u$ is even therefore $P(w)$ is also even. Moreover, integrating by parts we rewrite $P(w)$ as
\begin{equation*}
P(w)= \begin{cases}
\displaystyle e^{\frac{ix^2}{2}}-1-e^{\frac{ix^2}{2}}\int_0^x \frac{y e^{-\frac{iy^2}{2}}}{\Re(u(y))}\,dy,\quad x\in [0,1],\\
\displaystyle e^{\frac{ix^2}{2}}\left(1-\int_0^1 \frac{y e^{-\frac{iy^2}{2}}}{\Re(u(y))}\,dy\right)-1\\
\displaystyle -i\left(\frac{1}{\Re(u(x))}-\frac{e^{\frac{i(x^2-1)}{2}}}{\Re(u(1))}\right)-ie^{\frac{ix^2}{2}}\int_1^{x}
\frac{e^{-\frac{iy^2}{2}}\Re(u'(y))}{\Re(u(y))^2}\,dy,\quad x\in [1,+\infty).\end{cases}\end{equation*}

We first show that $P(E)\subset E$. 
 We set for $y\geq 0$
$$f(y)=\Re\left(\alpha+\int_{y}^{+\infty}\frac{w(z)}{z^2}\,dz\right)$$
so that $$\Re(u(y))=1+y f(y).$$
Then $f\in C^1(\R_+)$ and by \eqref{ineq:sup} it satisfies
\begin{equation}\label{bounds:g}|f'(y)|\leq \frac{\alpha}{4y^2},\quad |f(y)|\leq \frac{3\alpha}{2},\quad f(y)\geq \frac{\alpha}{2}>0.\end{equation}

On the one hand, {this yields $\Re(u(y))\geq 1$ so $P(w)(0)=P(w)'(0)=0$ and }
\begin{equation*}
\begin{split}
\int_0^1 \frac{y}{\Re(u(y))}\,dy\leq 1,
\end{split}
\end{equation*}
therefore
\begin{equation*}
\|P(w)\|_{L^\infty([0,1])}\leq 3.
\end{equation*}

On the other hand, using \eqref{bounds:g}  we get
\begin{equation*}
\begin{split}
\int_1^{+\infty} \frac{|\Re(u')|}{\Re(u)^2}\,dy&\leq \int_1^{+\infty} \frac{y|f'|+f}{(1+yf)^2}\,dy\leq \int_1^{+\infty} \left(\frac{y|f'|}{2y f}+\frac{f}{y^2f^2}\right)\,dy\\
&\leq \int_{1}^{+\infty}\frac{|f'|}{2f}\,dy+\frac{2}{\alpha}\int_1^{+\infty}\frac{dy}{y^2}\\
&\leq \frac{\alpha+8}{4\alpha}.
\end{split}
\end{equation*}
It follows that
\begin{equation*}
\begin{split}
\|P(w)\|_{L^\infty([1,+\infty))}\leq 5+\frac{\alpha+8}{4\alpha}.
\end{split}
\end{equation*}
Next, we have
\begin{equation*}\begin{split}
P'(w)=ix(1+P(w))-\frac{x}{\Re(u)}
\end{split}
\end{equation*}
so that combining the previous estimates
\begin{equation*}\begin{split}
\||x|^{-1}P'(w)\|_{L^\infty}\leq 7+ \frac{\alpha+8}{4\alpha}.
\end{split}
\end{equation*}

In conclusion, we obtain
\begin{equation*}
\|P(w)\|_{L^\infty}+\||x|^{-1}P'(w)\|_{L^\infty}\leq \frac{\alpha}{4}
\end{equation*}
provided that
\begin{equation}
\label{assumption-2}
12+\frac{\alpha+8}{2\alpha}\leq \frac{\alpha}{4},
\end{equation} which holds for $\alpha>K_0$ sufficiently large.

We then show that $P$ is a contraction on $E$. For $w_1,w_2\in E$ we set $f_i(y)=\Re\left(\alpha+\int_y^{+\infty} w_i/z^2\,dz\right)$, $i=1,2$. We have
\begin{equation*}
\begin{split}
P(w_1)(x)-P(w_2)(x)=-e^{\frac{ix^2}{2}}
\int_0^{x}ye^{-\frac{iy^2}{2}}g(y)\, dy,
\end{split}
\end{equation*}
where 
\begin{equation*}\begin{split}
g(y)&=\frac{{y}(f_2(y)-f_1(y))}{(1+{y}f_1(y))(1+{y}f_2(y))},\end{split}\end{equation*} so that
\begin{equation*}\begin{split}
g'(y)&=\frac{{y}(f_2'-f_1')+{f_2-f_1}}{(1+{y}f_1)(1+{y}f_2)}\\
&-{y}(f_2-f_1)\frac{(1+yf_1)(yf'_2+f_2)+(1+yf_1)(yf'_2+f_2)}
{(1+{y}f_1)^2(1+{y}f_2)^2}.
\end{split}
\end{equation*} 
Since $$f_1(y)-f_2(y)=\Re\left(\int_y^{+\infty} \frac{w_1-w_2}{z^2}\,dz\right)$$ we obtain by \eqref{ineq:int} for $y>0$
\begin{equation}\begin{split}
\label{bounds:diff}
|f_1(y)-f_2(y)|\leq d(w_1,w_2),\quad |f'_1(y)-f'_2(y)|\leq \frac{d(w_1,w_2)}{y^2},
\end{split}
\end{equation}
while
\begin{equation}
\label{bounds:diff-2}\frac{\alpha}{2}\leq f_i(y)\leq \frac{3\alpha}{2},\quad i=1,2.
\end{equation}
Using the inequality $1+yf_i(y)\geq yf_i(y)\geq y\alpha/2$, we obtain by \eqref{bounds:diff} and \eqref{bounds:diff-2} 
\begin{equation}\label{h1}
\begin{split}
|g(y)|&\leq \frac{y|f_2(y)-f_1(y)|}{yf_2(y)}\leq \frac{2}{\alpha}d(w_1,w_2),
\end{split}
\end{equation}
and there exists some $K_1>1$ such that
\begin{equation}\label{h2}
\begin{split}
|g'(y)|&\leq \frac{K_1}{\alpha^2y^{2}}d(w_1,w_2),\quad y \geq 1.
\end{split}
\end{equation}
Hence 
\begin{equation*}
\|P(w_1)-P(w_2)\|_{L^\infty([0,1])}\leq  \frac{1}{\alpha} d(w_1,w_2).
\end{equation*}
And for $x\geq 1$ an integration by parts yields
\begin{equation*}
\begin{split}
P(w_1)(x)-P(w_2)(x)&=-e^{\frac{ix^2}{2}}
\int_0^{1}ye^{-\frac{iy^2}{2}}hg(y)\, dy-ig(x)+ie^{\frac{i(x^2-1)}{2}}g(1)\\
&+ie^{\frac{ix^2}{2}}\int_1^{x} e^{-\frac{iy^2}{2}}g'(y)\,dy
\end{split}
\end{equation*}
hence in view of \eqref{h1}-\eqref{h2}
\begin{equation*}
\begin{split}
\|P(w_1)-P(w_2)\|_{L^\infty([1,+\infty))}&\leq \Big(\frac{5}{\alpha}+\frac{K_1}{\alpha^2} \Big)d(w_1,w_2).
\end{split}
\end{equation*}

On the other hand,
\begin{equation*}
P'(w_1)(x)-P'(w_2)(x)=ix\left(P(w_1)(x)-P(w_2)(x)\right)-xg(x)\end{equation*}
therefore using again \eqref{bounds:diff}-\eqref{bounds:diff-2} we get
\begin{equation*}
\| x^{-1}(P'(w_1)-P'(w_2))\|_{L^\infty}\leq \left(\frac{7}{\alpha}+\frac{K_1}{\alpha^2}\right)\;d(w_1,w_2).
\end{equation*}

Finally, 
\begin{equation*}
d(P(w_1),P(w_2))
\leq \left(\frac{12}{\alpha}+\frac{2K_1}{\alpha^2}\right)\:d(w_1,w_2),
\end{equation*}
which establishes that $P$ is a contraction provided that $\alpha$ satisfies
\begin{equation}\label{assumption-3}
\frac{12}{\alpha}+\frac{2 K_1}{\alpha^2}<1,
\end{equation} and increasing possibly $K_0$ we are ensured that both conditions \eqref{assumption-2} and \eqref{assumption-3} are satisfied. This concludes the proof of Theorem \ref{thm:ak}.

\end{proof}

\section{Numerical simulations\label{secnum}}

To compute an approximation of the system \eqref{KMD2gen}, we use a splitting algorithm between the nonlinear and linear parts. Indeed, the solution of the linear part 
\begin{equation}
\label{eq:libre}
\begin{cases}
\displaystyle i\partial_t \Psi_1+\alpha_1\partial_\sigma^2\Psi_1=0,\\
\displaystyle i\partial_t \Psi_2-\alpha_2\partial_\sigma^2\Psi_2=0,\\
\end{cases}
\end{equation}
is given explicitly in Fourier, and can be very easily computed using the Fast Fourier Transform algorithm. On the other hand, the solution of the nonlinear part 
\begin{equation*}
\begin{cases}
\displaystyle i\partial_t \Psi_1-\frac{\Psi_1-\Psi_2}{|\Psi_1-\Psi_2|^2}=0,\\
\displaystyle i\partial_t \Psi_2-\frac{\Psi_1-\Psi_2}{|\Psi_1-\Psi_2|^2}=0,\\
\end{cases}
\end{equation*}
is given explicitly by the formula
$$
 \Psi_k(t,x) = \Psi_k(0,x) - i t \frac{\Psi_1(0,x)-\Psi_2(0,x)}{|\Psi_1(0,x)-\Psi_2(0,x)|^2}.
$$
Denoting by $\varphi_k^t(\Psi_1(0),\Psi_2(0))$ this application, 
the algorithm used to compute the solutions are thus based on the classical Lie approximation 
$$
\Psi_k(t) \simeq e^{i t \alpha_k \partial_\sigma^2} \circ \varphi_k^t(\Psi_1(0), \Psi_2(0)), 
$$
as well as the symmetric version (Strang splitting) 
$$
\Psi_k(t) \simeq \varphi_k^{t/2} \circ e^{i t  \alpha_k \partial_\sigma^2} \circ \varphi_k^{t/2} (\Psi_1(0), \Psi_2(0)), 
$$
To discretize in space these formulae, we take a large periodic box $L[-\pi,\pi]$ using equidistant grid points $x_k = L k \pi/K$, for $k \in B_K := \{K/2 - 1,K/2\}$. The number $K$ is the number of nodes. Then we can evaluate the flow $\varphi^t$ at the fix nodes $x_k$, while the approximation of the free equation $e^{i t \alpha_k}$ is made using the discrete Fourier transformation: 
$$
\forall\, k \in B_K, \quad v_k = (\mathcal{F}_K)_k := \frac{1}{K}\sum_{\ell \in B_K} u_\ell e^{- i k x_\ell}
$$
and the calculation $v_k(t) = e^{-i t k^2L^{-2}} v_k(0)$ as solution of \eqref{eq:libre} in Fourier variables. We then go back to the $x$-variable using the inverse of $\mathcal{F}_K$ that, together with $\mathcal{F}_K$, can be computed in $K \log K$ operations. 

The algorithms obtained are then symplectic, and used with a time stepsize $\tau$ satisfying the CFL condition $\tau K^2L^{-2} < \pi$ which is known to be a necessary condition to ensure the existence of a modified energy as well as some stability results in simpler situations (cubic Schr\"odinger equation for instance, see \cite{BFG,F11}).  

\subsection{{Opposite core parameters}}

\begin{figure}
   \includegraphics[height=4.2cm]{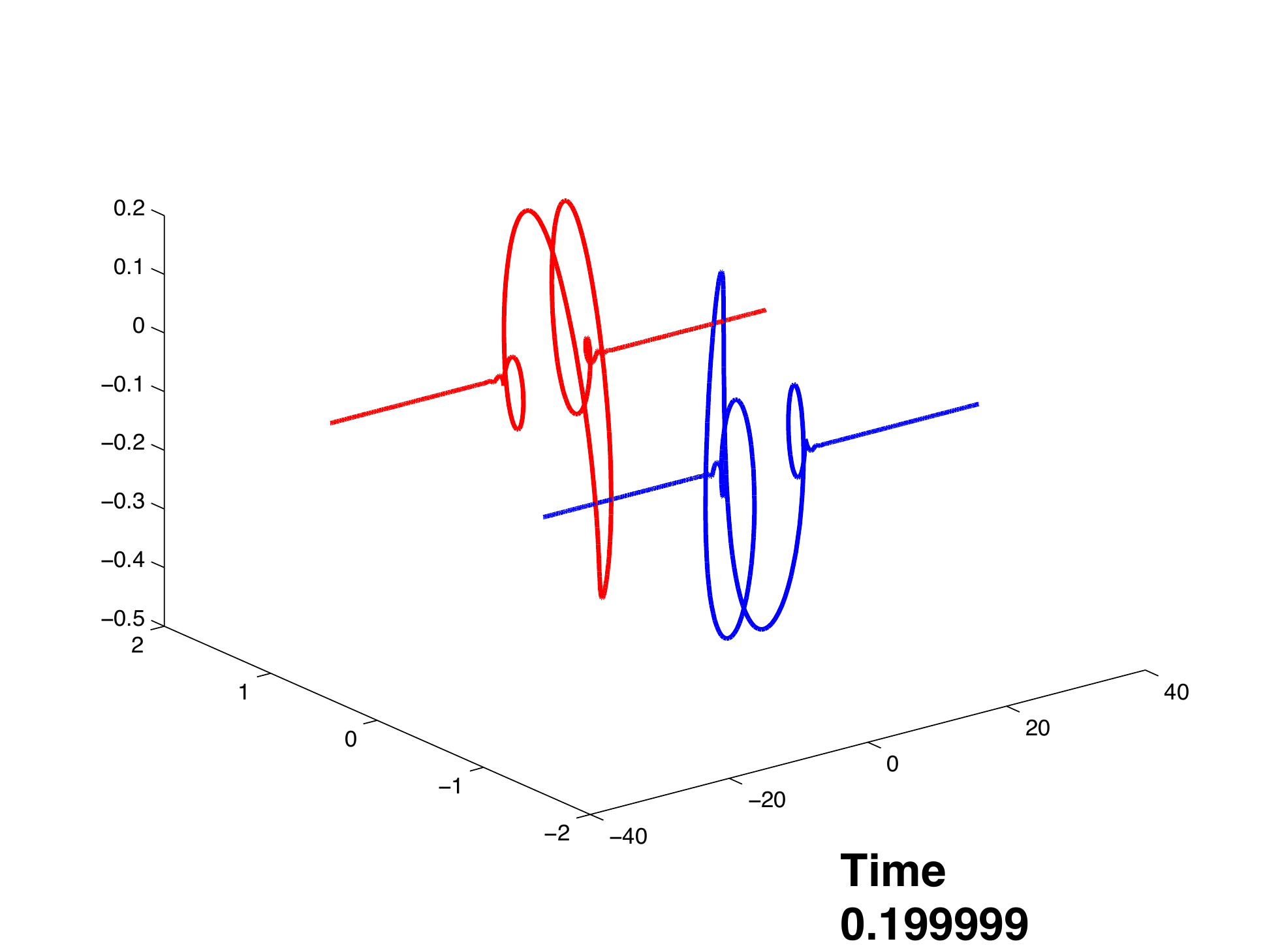}
   \includegraphics[height=4.2cm]{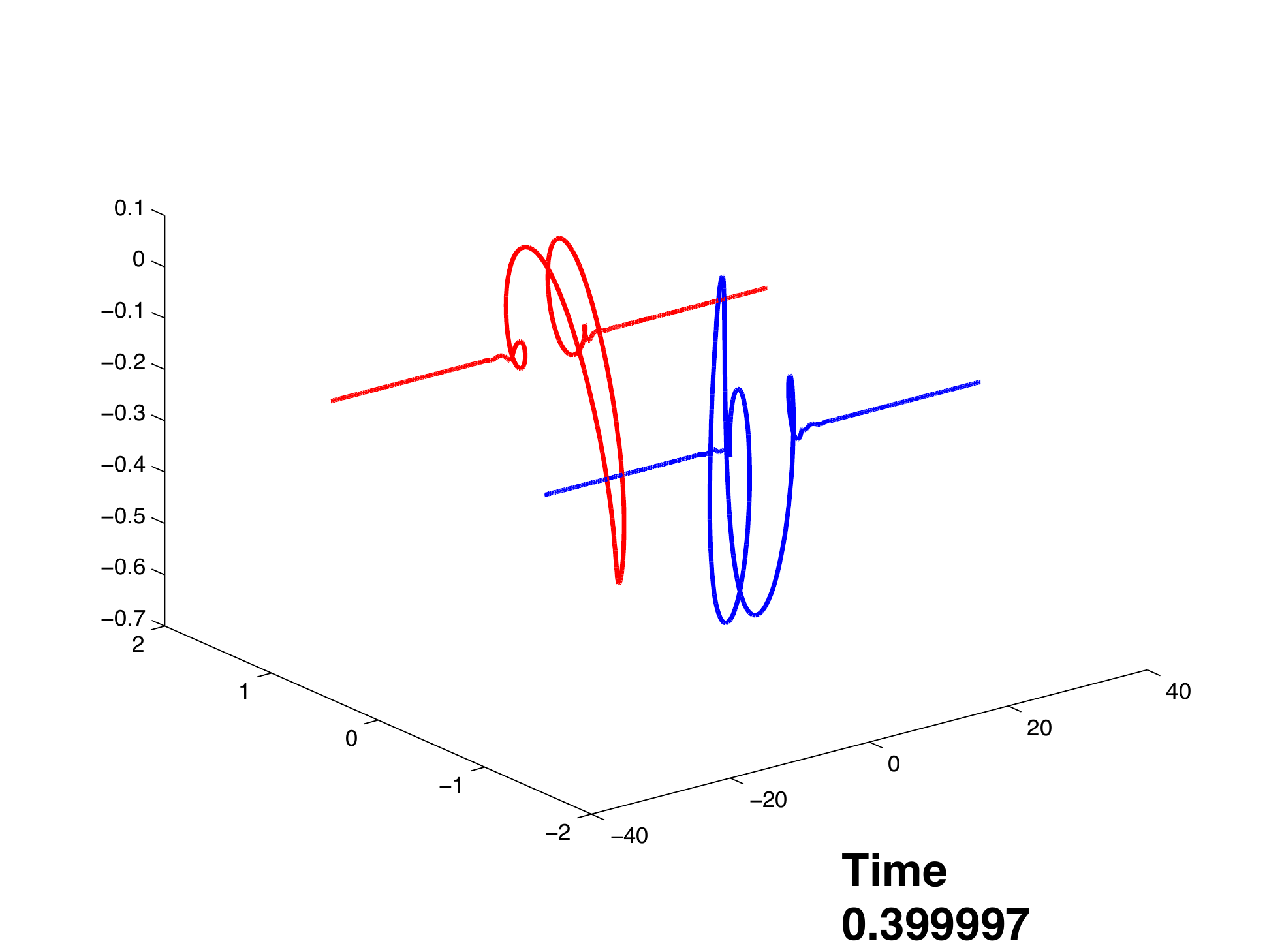} 
   \includegraphics[height=4.2cm]{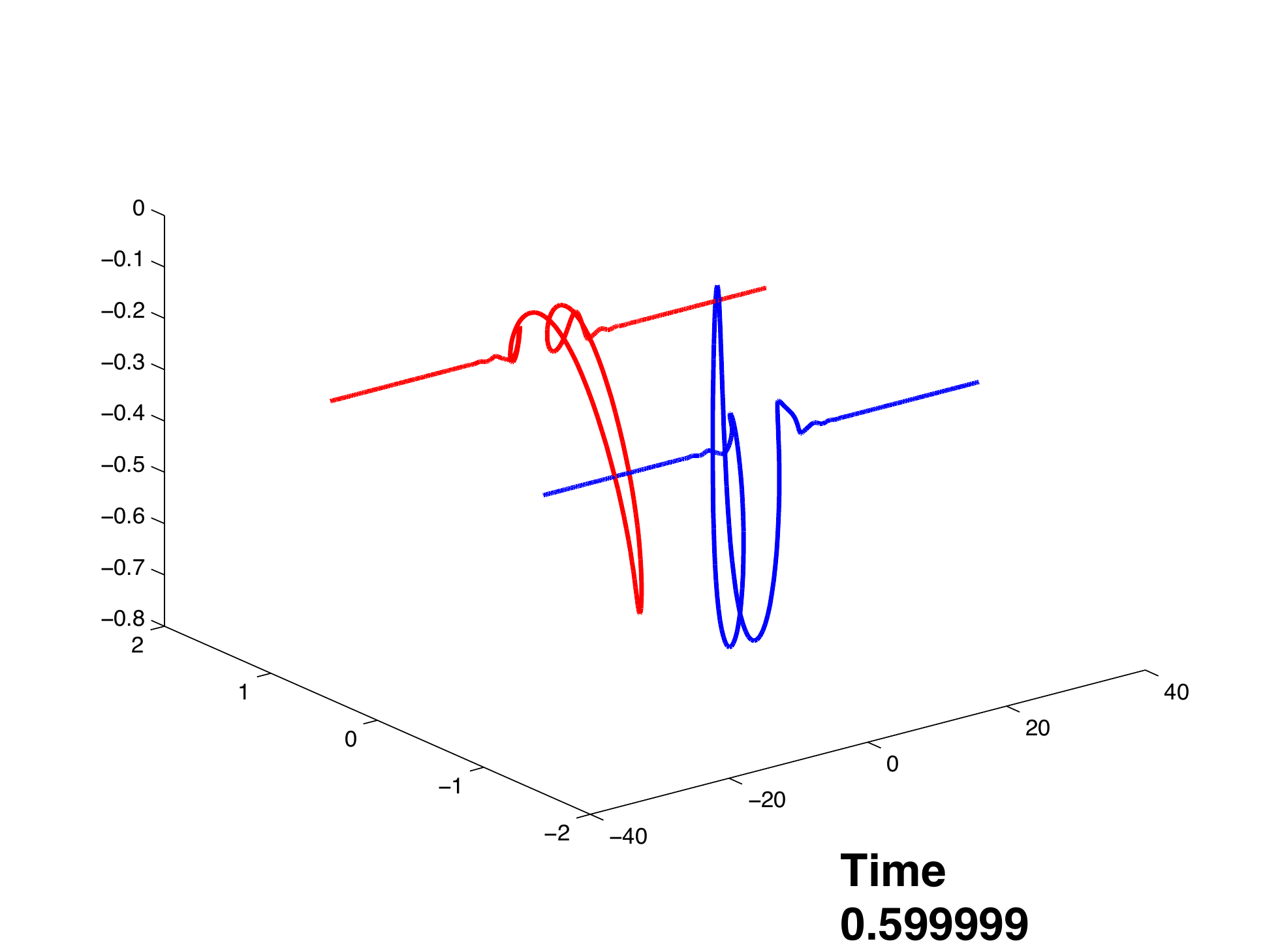},
   \includegraphics[height=4.2cm]{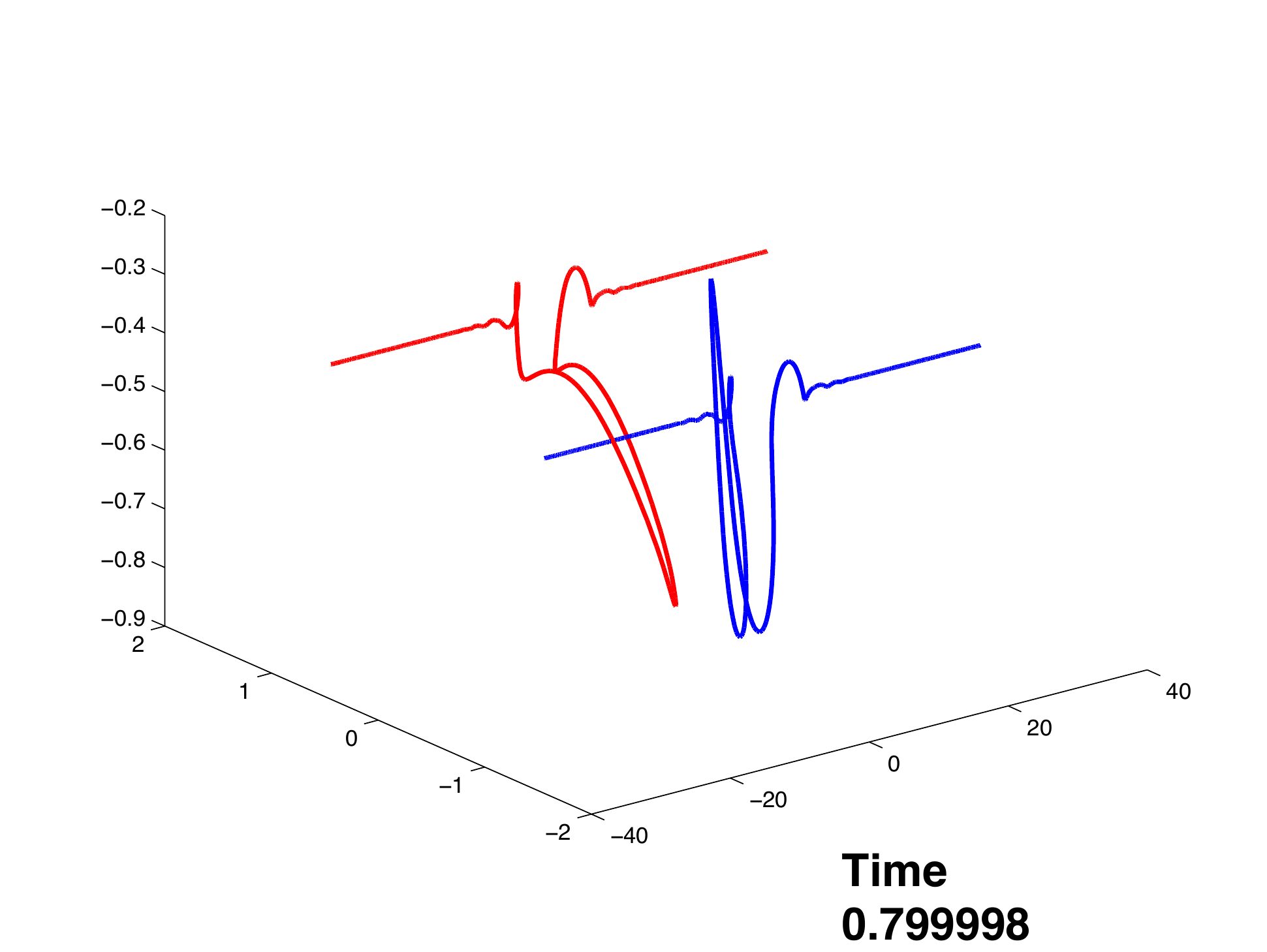}
   \includegraphics[height=4.2cm]{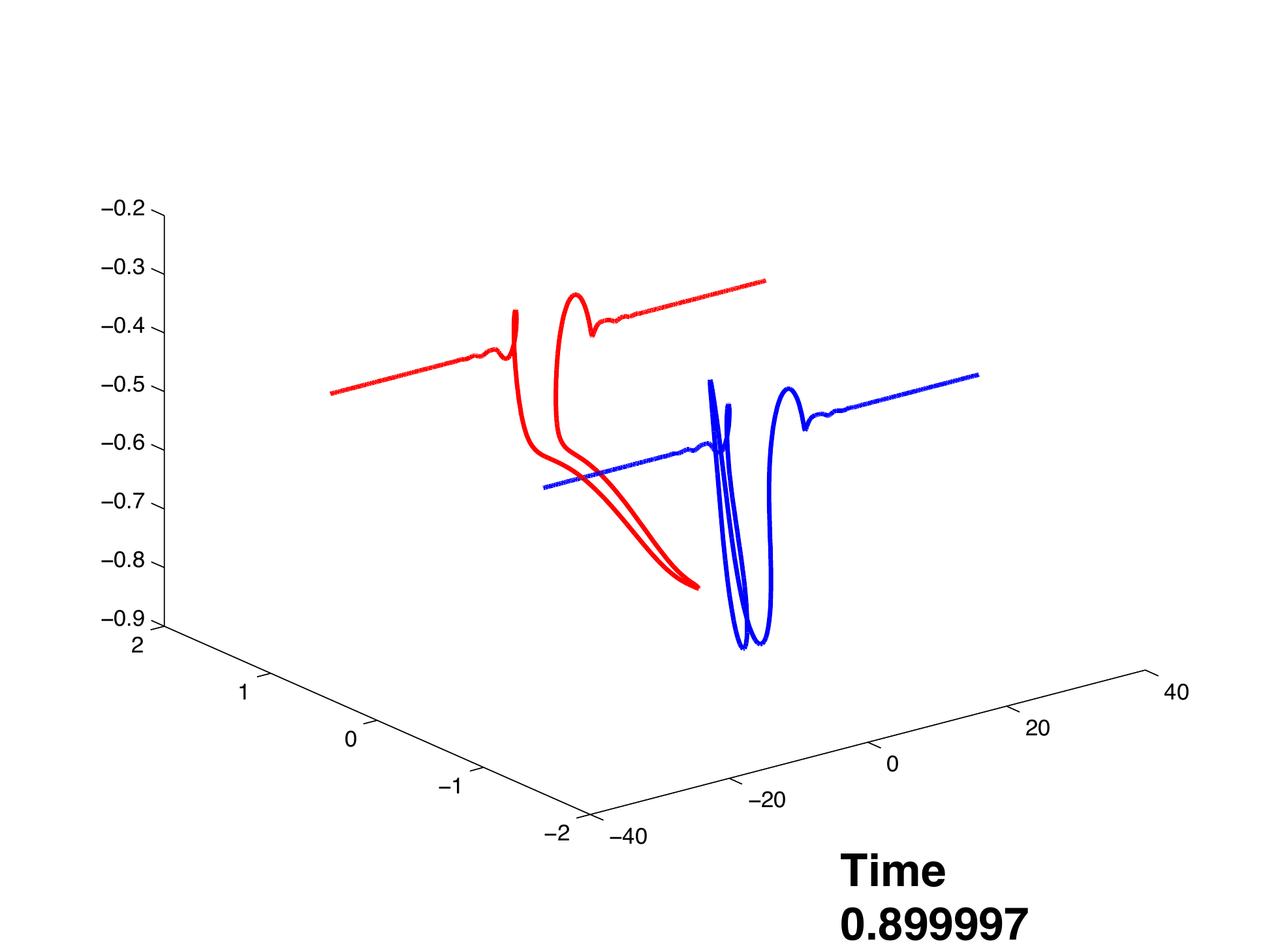}
   \includegraphics[height=4.2cm]{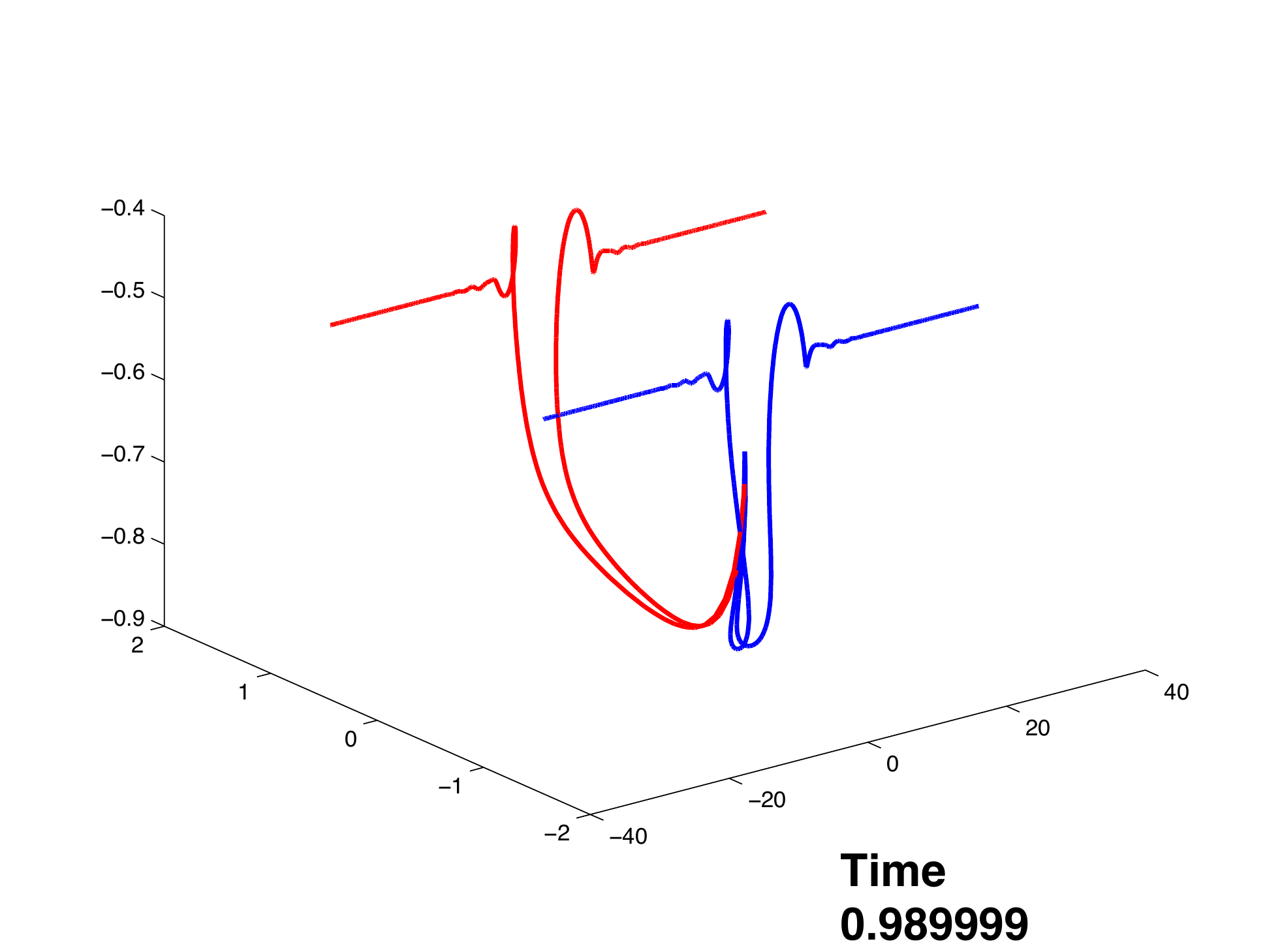}
  \caption{Collision in the case $\alpha_1=-\alpha_2=1$, $\Psi_1(0)=1-G$, $\Psi_2(0)=-1+G(\sigma)$, {no symmetry  kept by $(\Psi_1,\Psi_2)(t,\sigma)$}.}
    \label{opposite}
\end{figure}

We consider now the equation {\eqref{KMD2gen},} in the case where $\Gamma_1 = -\Gamma_2 = 1$ and $\alpha_1 = -\alpha_2 = 1$. We take $K = 1024$, and the stepsize $\tau = \pi K^{-2}$. We take $L = 10$ (so that the CFL number is of order $10^{-2}$), and the initial data $\Psi_1(0) = 1 - G$ and $\Psi_2(0) = -1 + G$. The {non-symmetric} evolution of the corresponding filaments is depicted in Figure \ref{opposite}. We observe that the collision occurs at time $t = 1$ as predicted by Proposition \ref{prop21}. 

\subsection{{Positive core parameters}}

\begin{figure}
   \includegraphics[height=4.2cm]{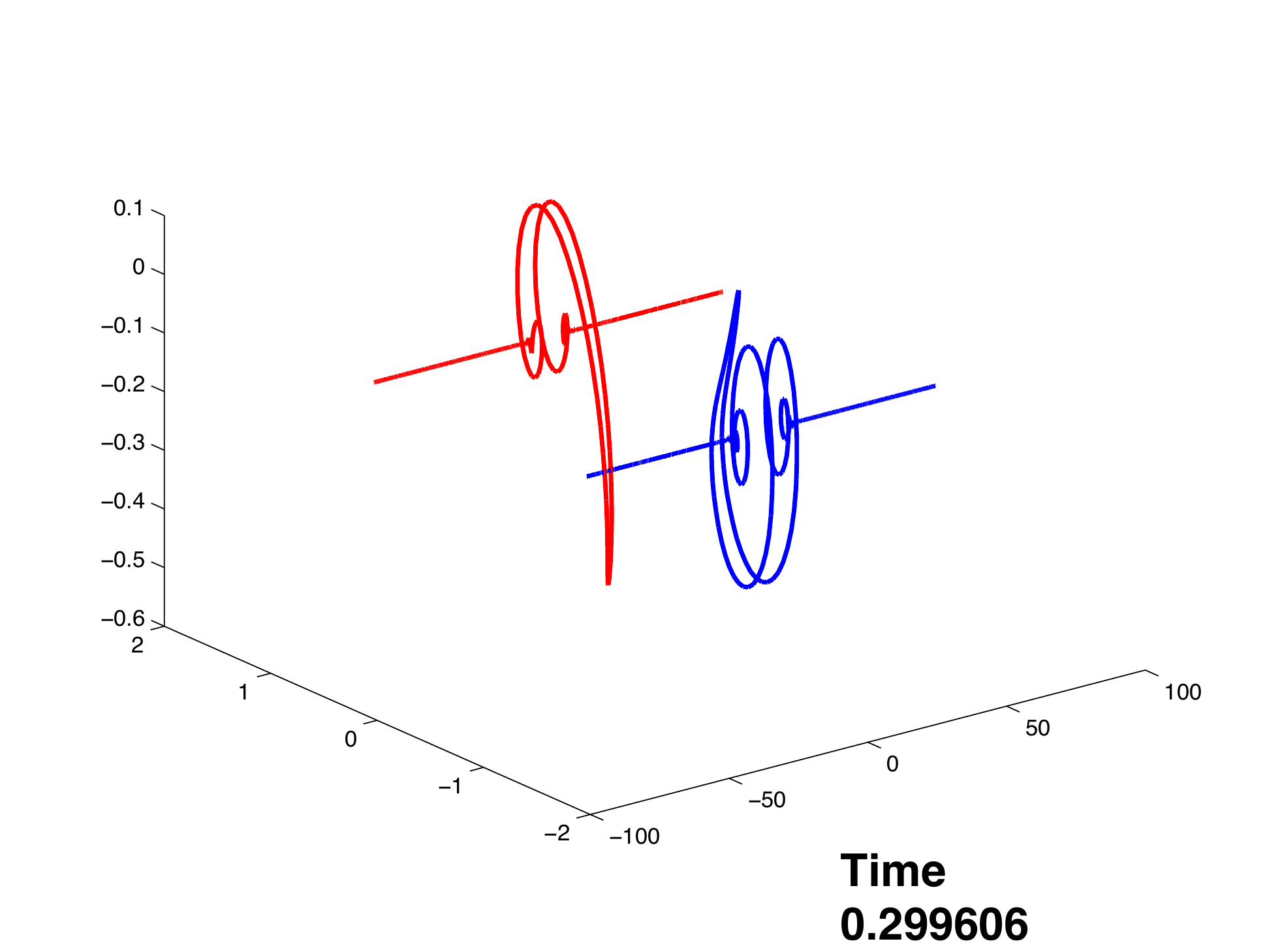}
   \includegraphics[height=4.2cm]{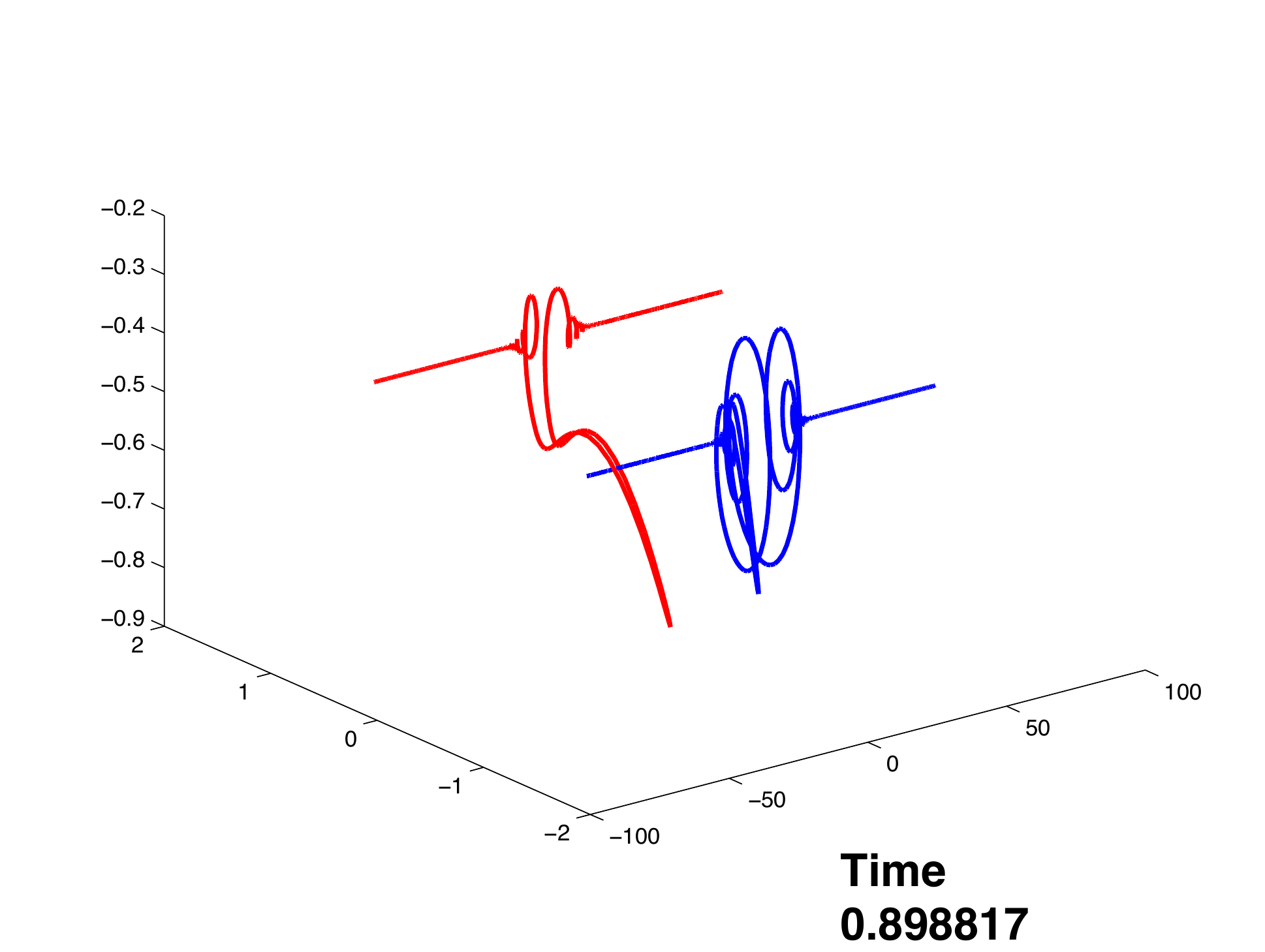} 
   \includegraphics[height=4.2cm]{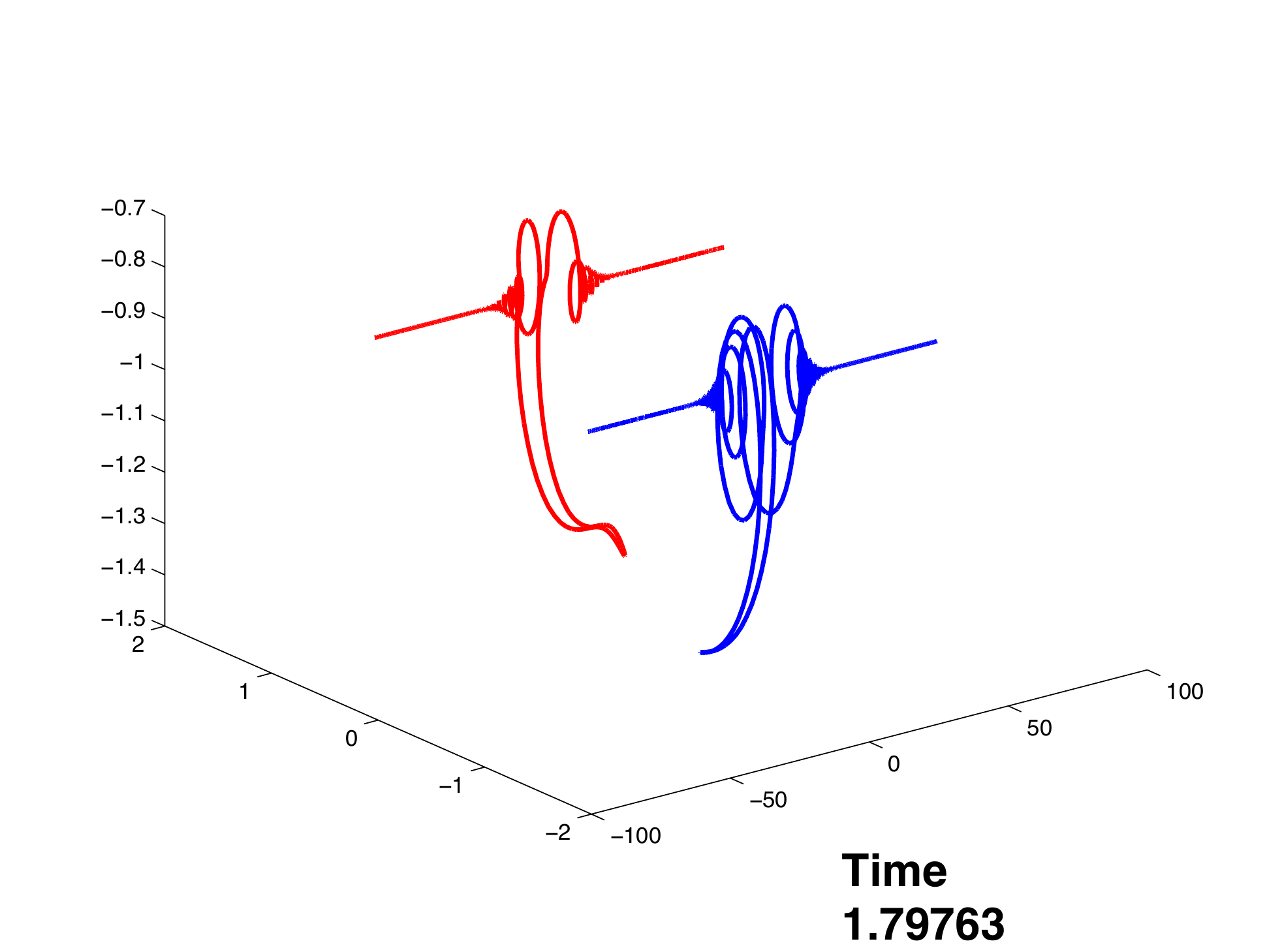},
   \includegraphics[height=4.2cm]{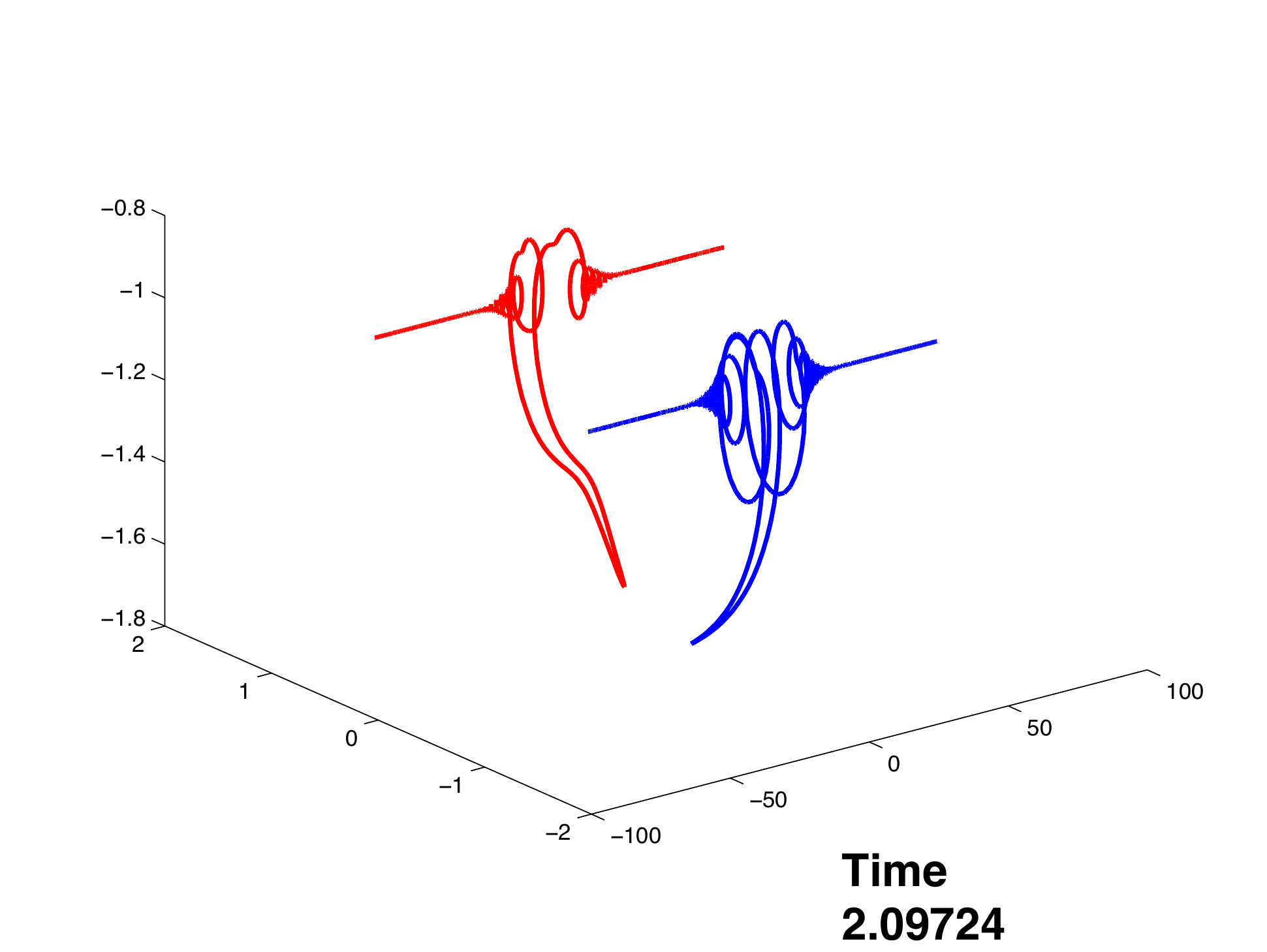}
   \includegraphics[height=4.2cm]{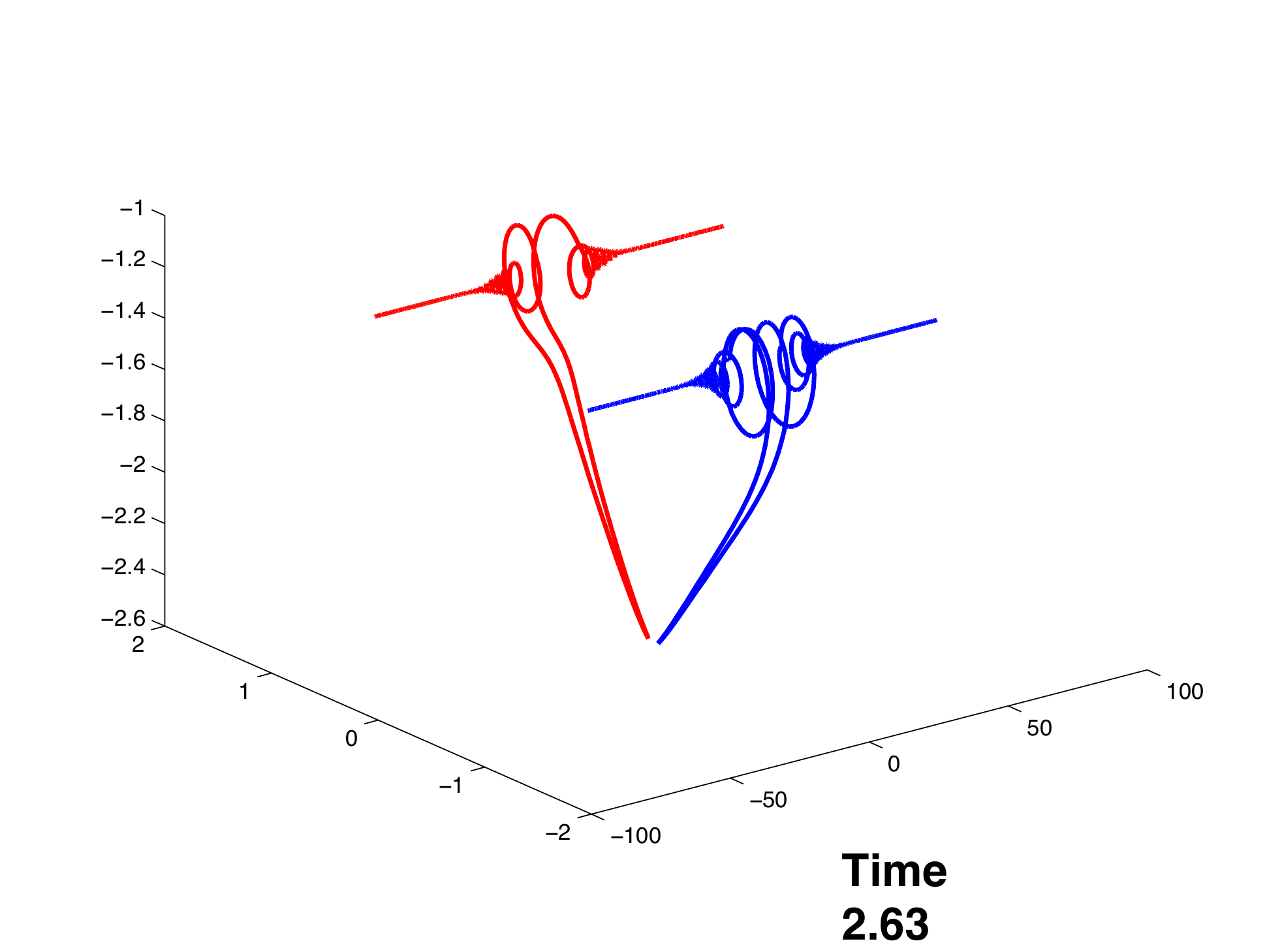}
   \includegraphics[height=4.2cm]{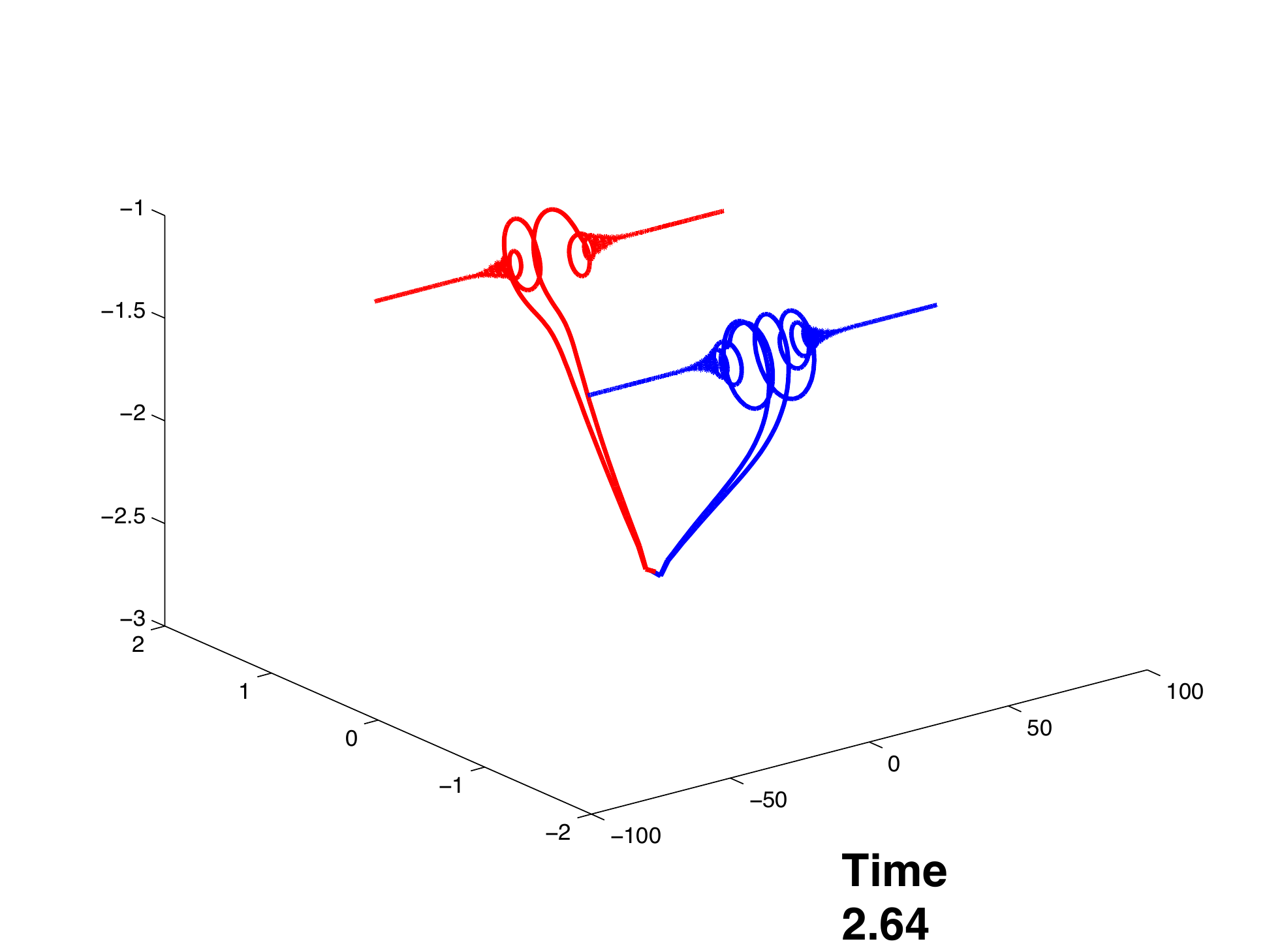}
  \caption{Collision in the case $\alpha_1=\alpha_2=1$, $\Psi_1(0)=1-G$, $\Psi_2(0)=-1+G$, no symmetry kept by  $(\Psi_1,\Psi_2)(t)$.}
    \label{samenosym}
\end{figure}

We now consider exactly the same initial condition as before, but we just change the sign of $\alpha_2$ setting $\alpha_1=\alpha_2=1$, so the evolution is governed by \eqref{KMD2same}. As mentioned above, the solution starting with the initial data $\Psi_1(0,\sigma) = 1 - G(\sigma)$ and $\Psi_2(0,\sigma) = -1 + G(\sigma)$ leads to a collision behavior in $\sigma = 0$ at a time of order $t \simeq 2.64$, see Figure \ref{samenosym}. Note that the solution is more oscillatory compared to the situation where $\alpha_1{=-\alpha_2,}$ and the support slightly larger. Hence we have taken $L = 20$ instead of $L = 10$ in this simulation. 
  \\

\begin{figure}
  \includegraphics[height=4.2cm]{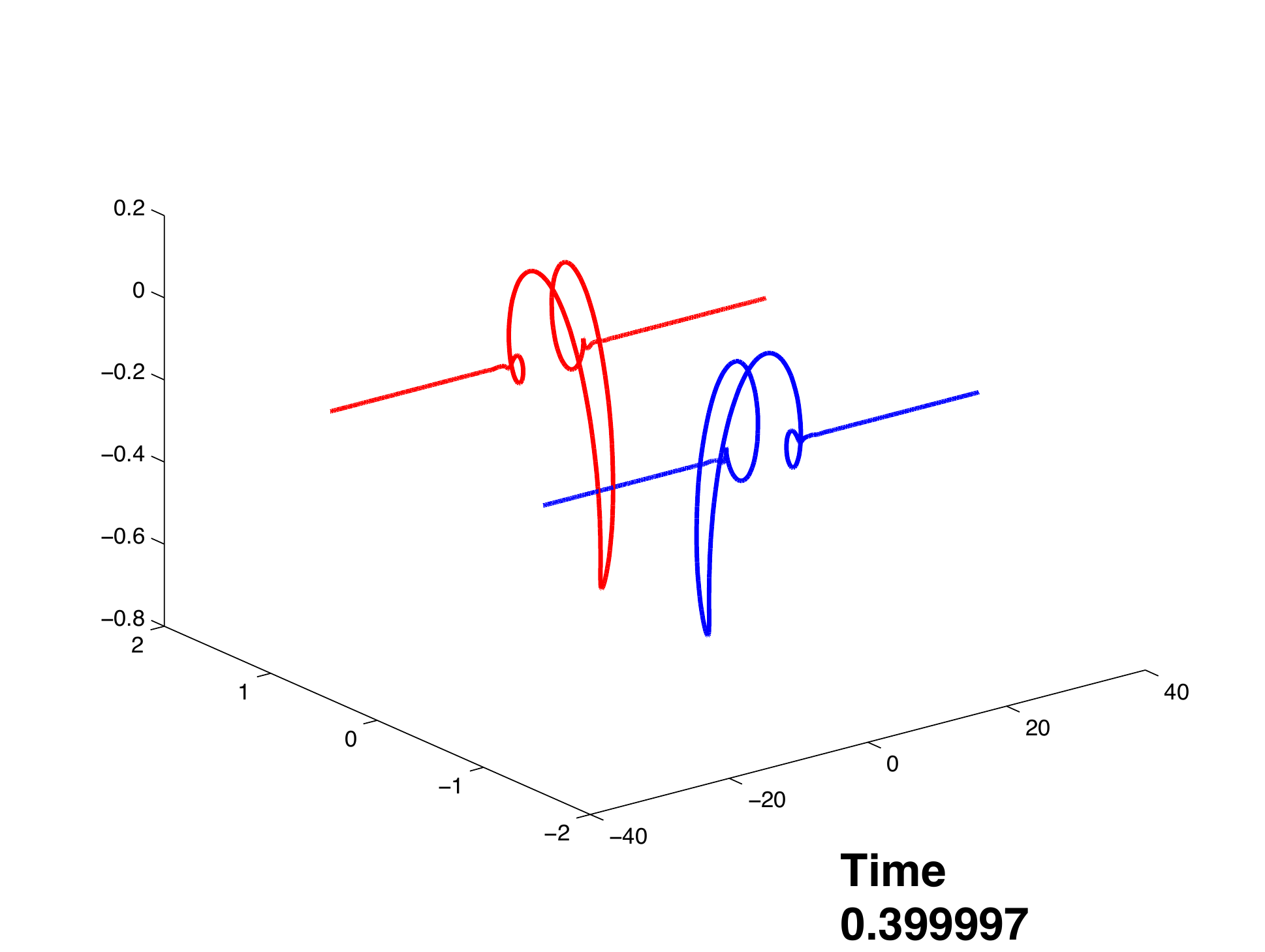}
   \includegraphics[height=4.2cm]{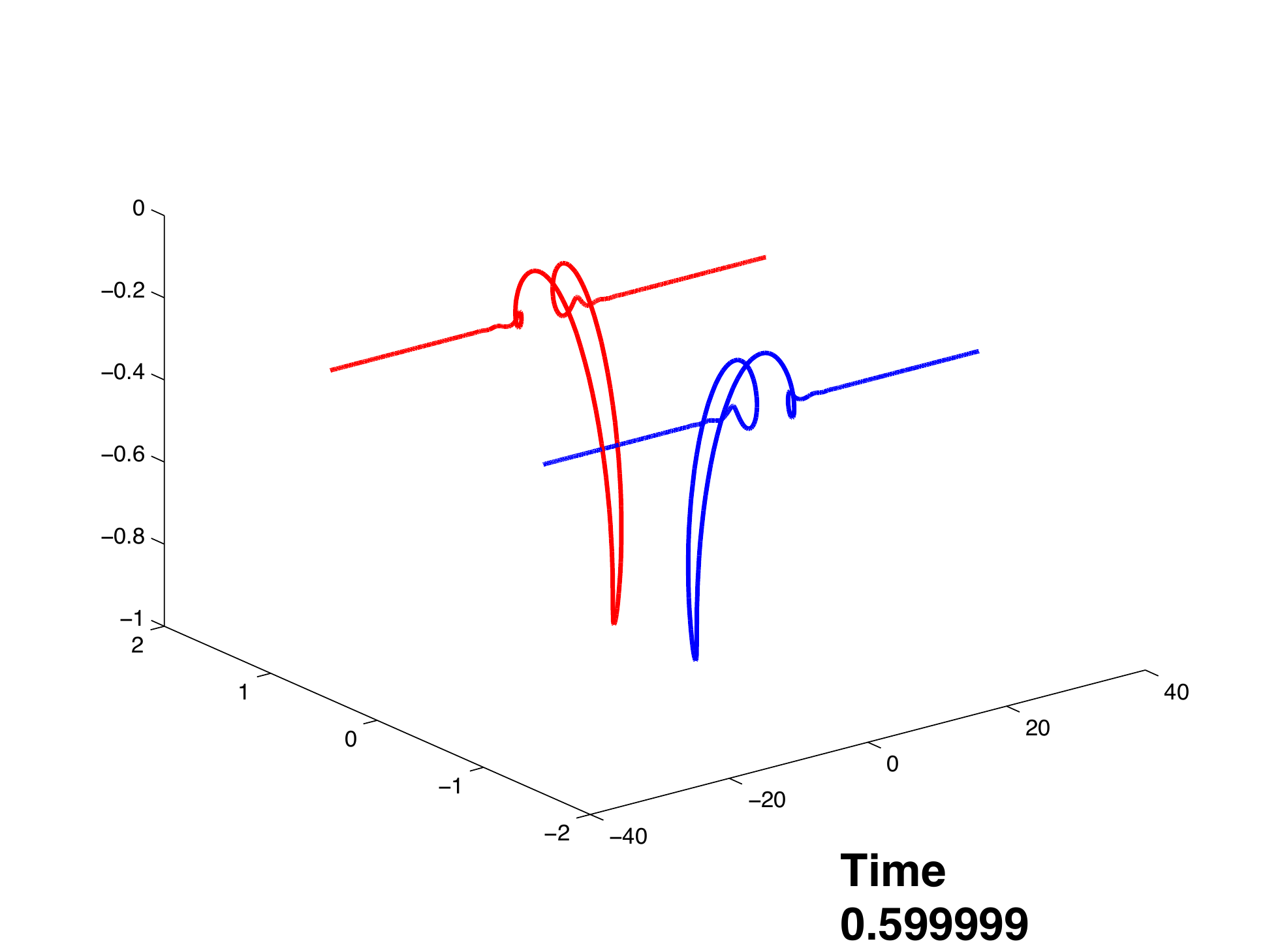} 
   \includegraphics[height=4.2cm]{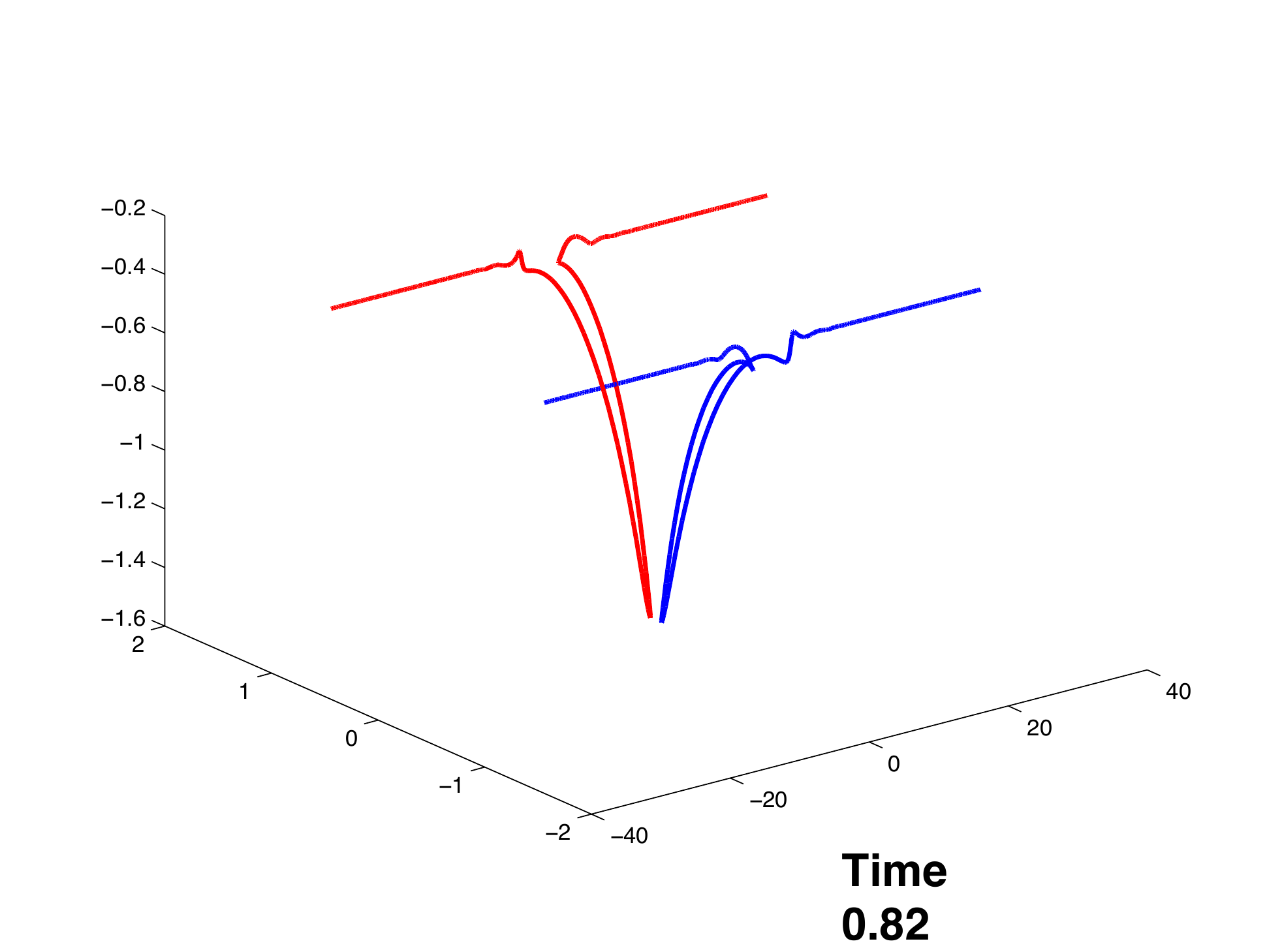},
   \includegraphics[height=4.2cm]{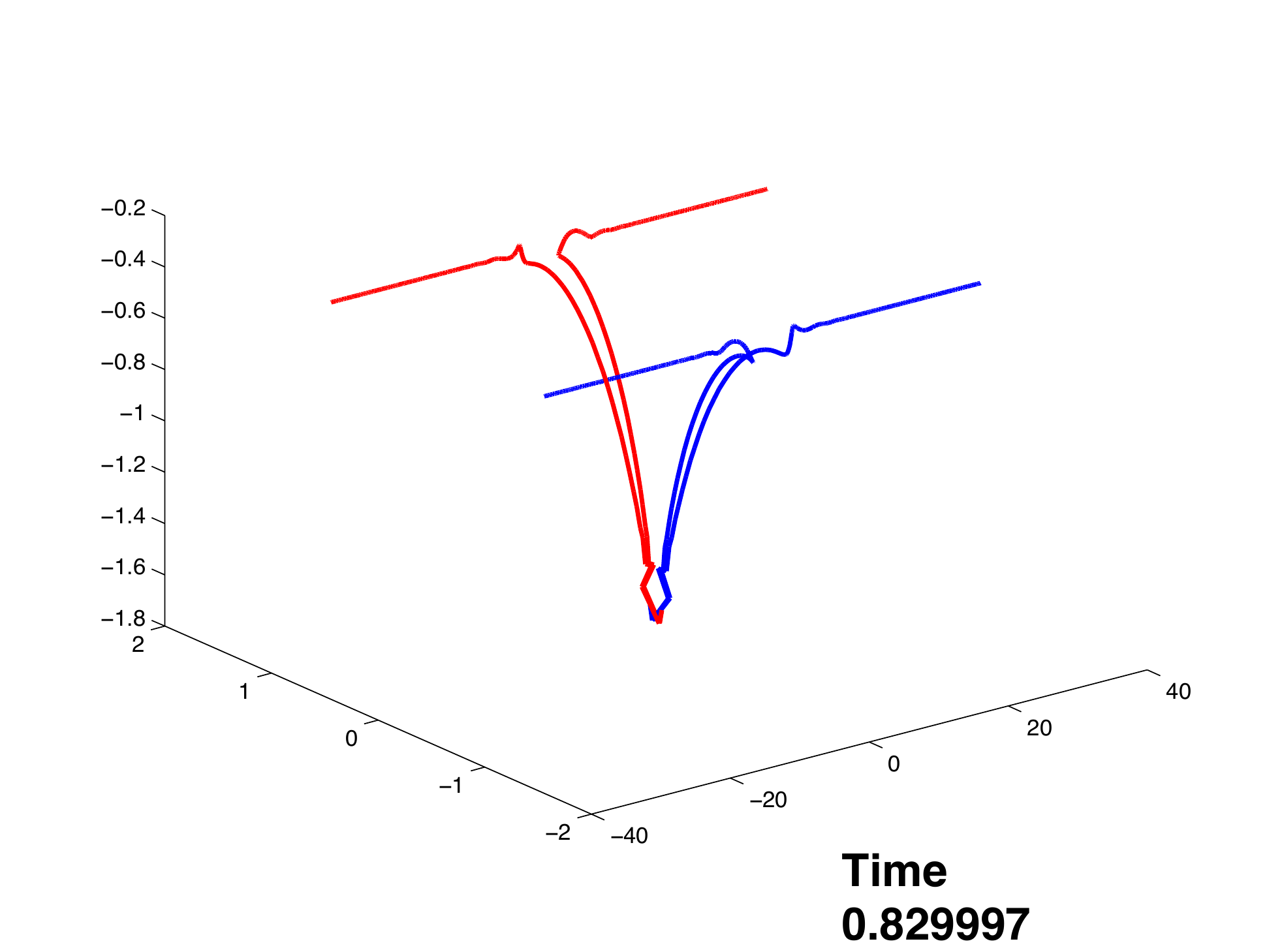}
  \caption{Collision in the case $\alpha_1=\alpha_2=1$, $\Psi_1(0)=1-G$, $\Psi_2(0)=-1+\overline{G}$, symmetry $\Psi_1(t)=-\overline{\Psi}_2(t)$.}
    \label{same}
\end{figure}

In a second simulation, we take as initial data $\Psi_1(0) = 1 - G$ and $\Psi_2(0) = -1 + \overline{G}=-\overline{\Psi}_1(0)$ so that $\Psi_2(t)=-\overline{\Psi}_1(t)$. With the same parameters as before (and with $L = 10$) we observe a {symmetric} collision at time $t \simeq 0.83$, see Figure \ref{same}. 
  \medskip

\begin{figure}
  \includegraphics[height=4.2cm]{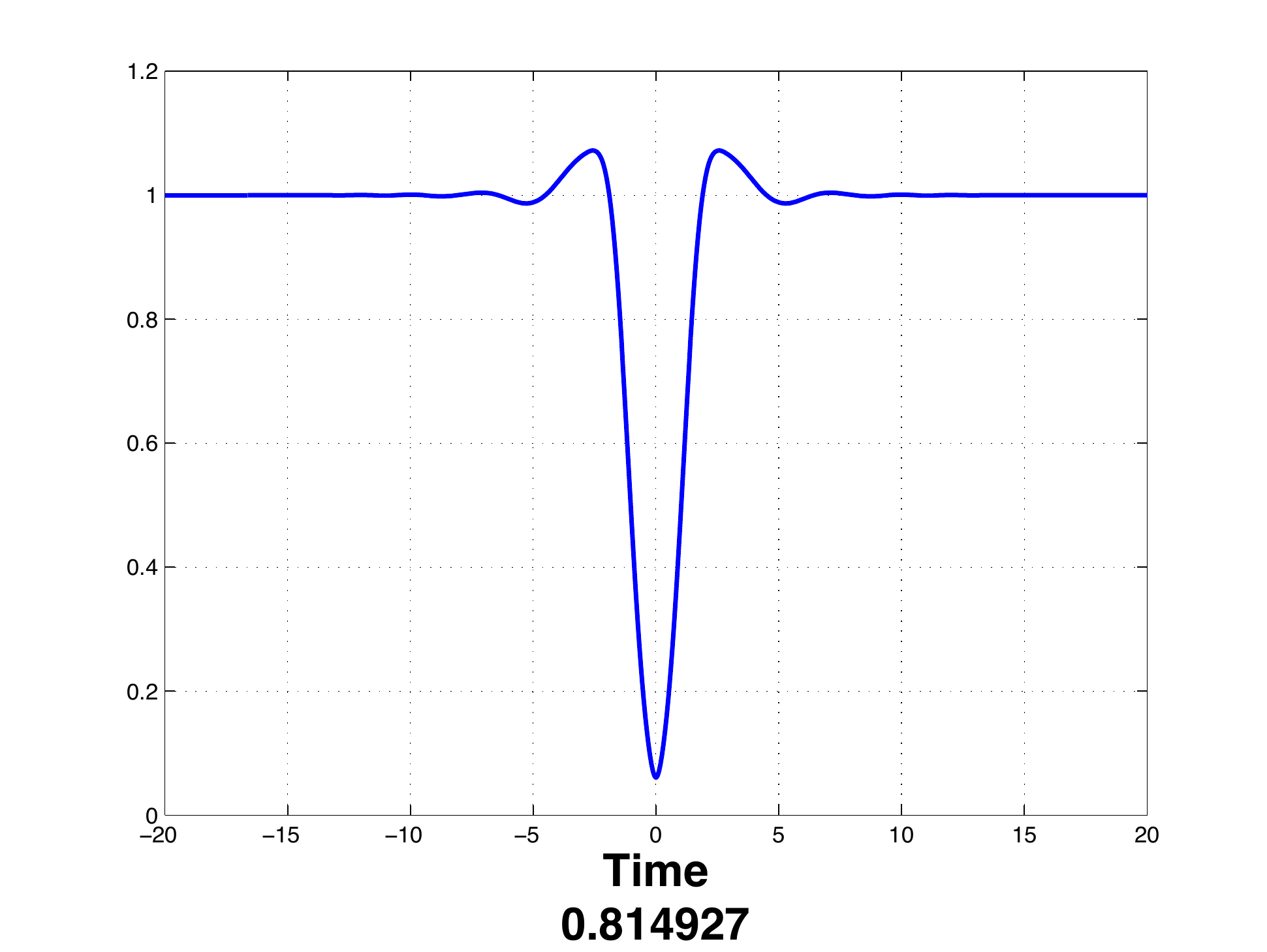}
   \includegraphics[height=4.2cm]{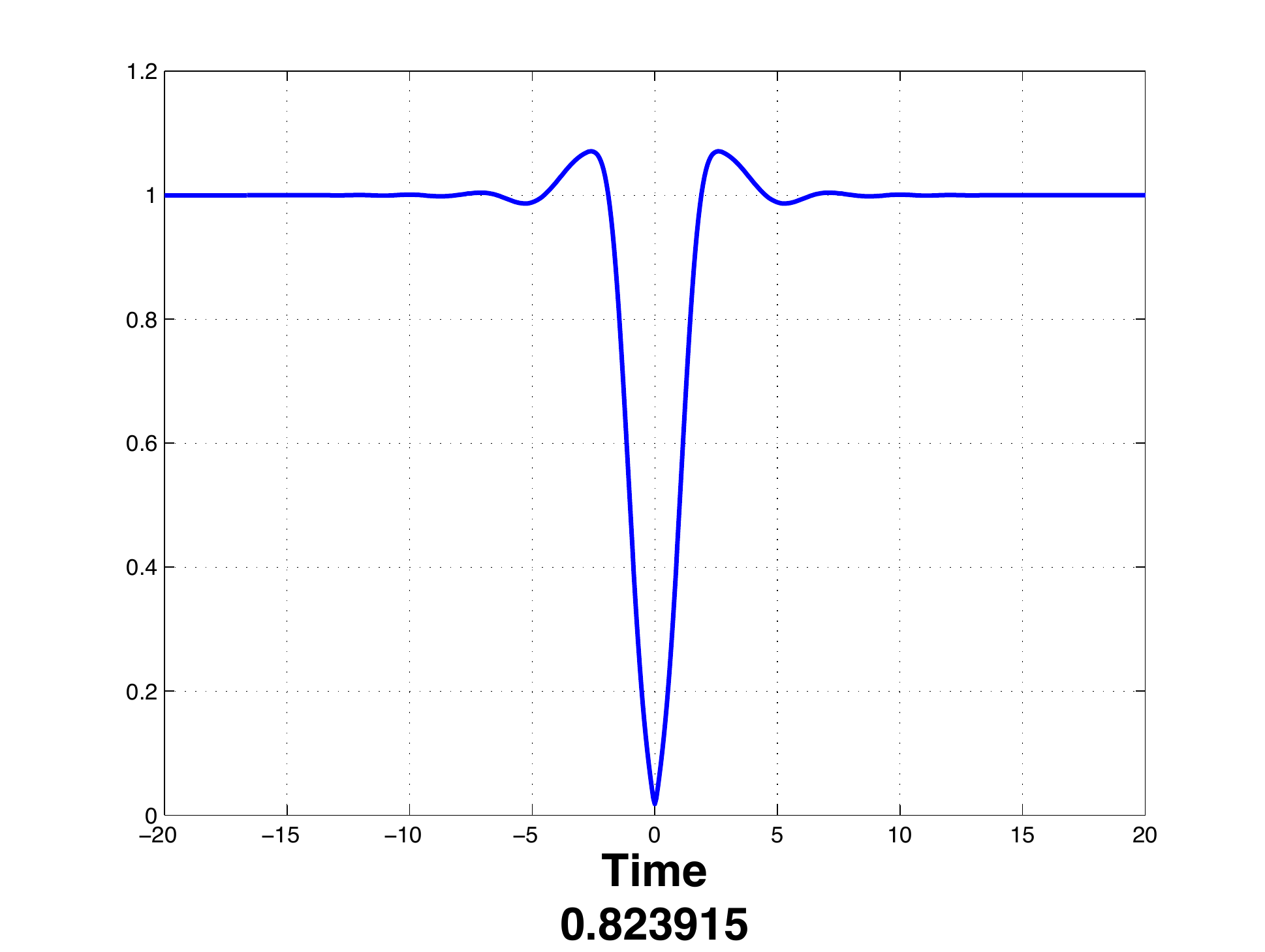} 
   \includegraphics[height=4.2cm]{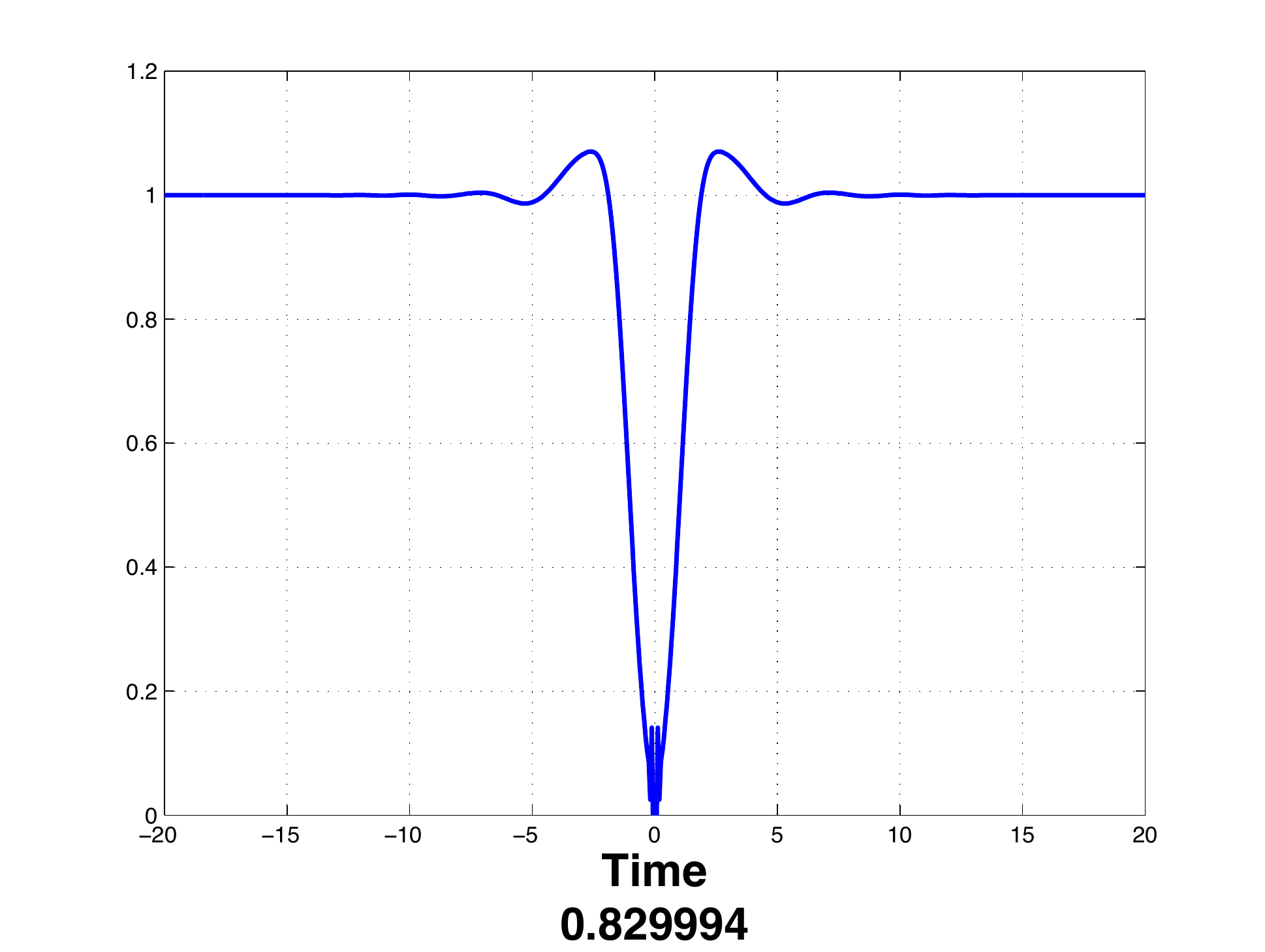},
    \includegraphics[height=4.2cm]{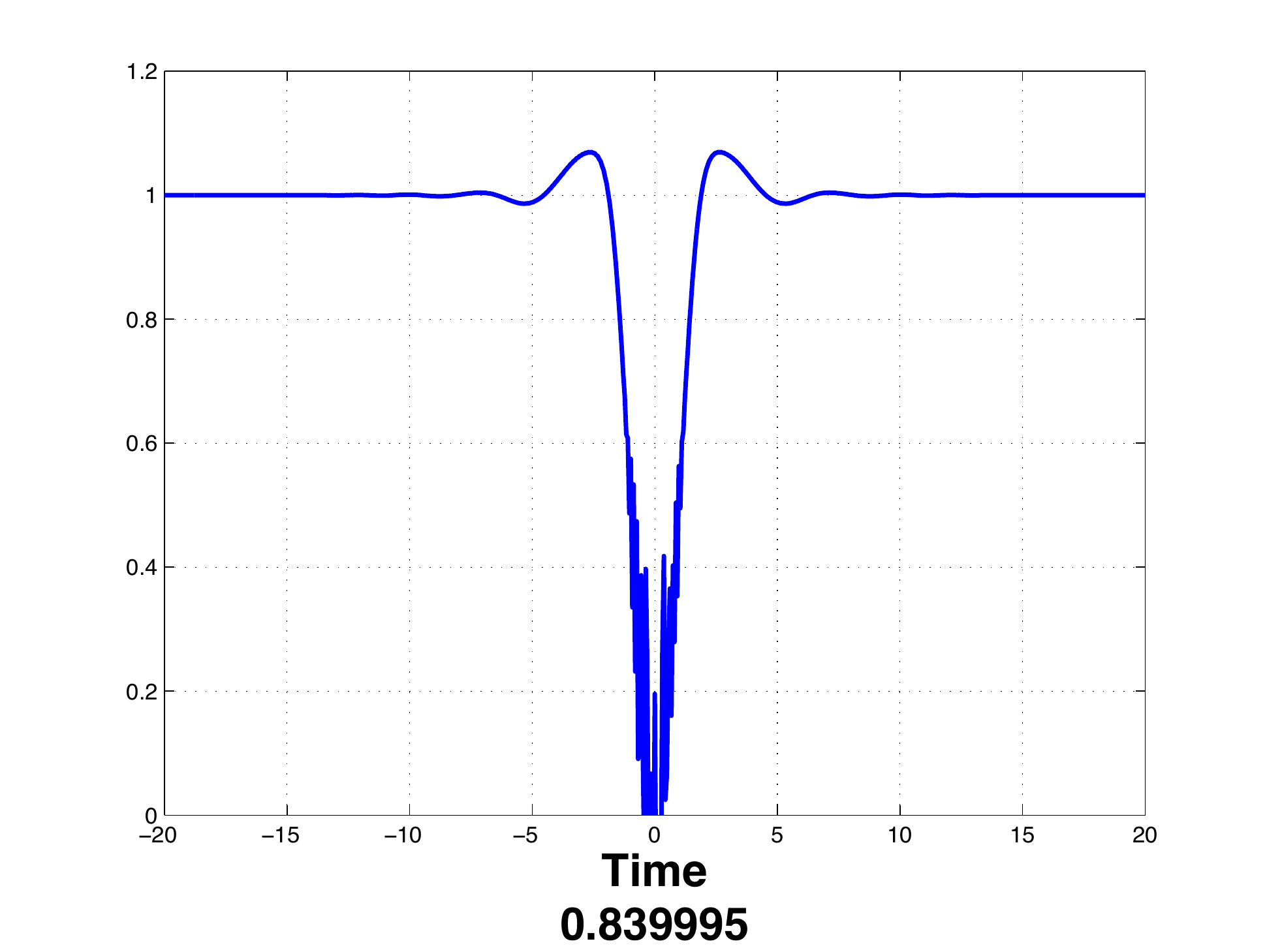}
  
  \caption{Real part of $\psi_1(t)$ for collision in the case  $\alpha_1=\alpha_2=1$, $\Psi_1(0)=1-G$, $\Psi_2(0)=-1+\overline{G}$, symmetry $\Psi_1(t)=-\overline{\Psi}_2(t)$.}
    \label{re1}
\end{figure}

To investigate the type of collision appearing in the symmetric situation $\Psi_1(0) ={-} \overline{\Psi}_2(0)$ we consider numerical approximations of Equation \eqref{eqpsi}. We first consider {again} the case where $\Psi_1(0) = 1 - G$. In Figure \ref{re1} we plot the evolution of the real part of the solution of \eqref{eqpsi} using the same splitting method (which extends straightforwardly to the this situation) and using again 1024 grid points and the same CFL condition (and $L = 10$). We observe that the collision (which means that the real part vanishes) occurs around the time $t = 0.83$, which confirms the previous numerical experiment. 

\begin{figure}
  \includegraphics[height=4.2cm]{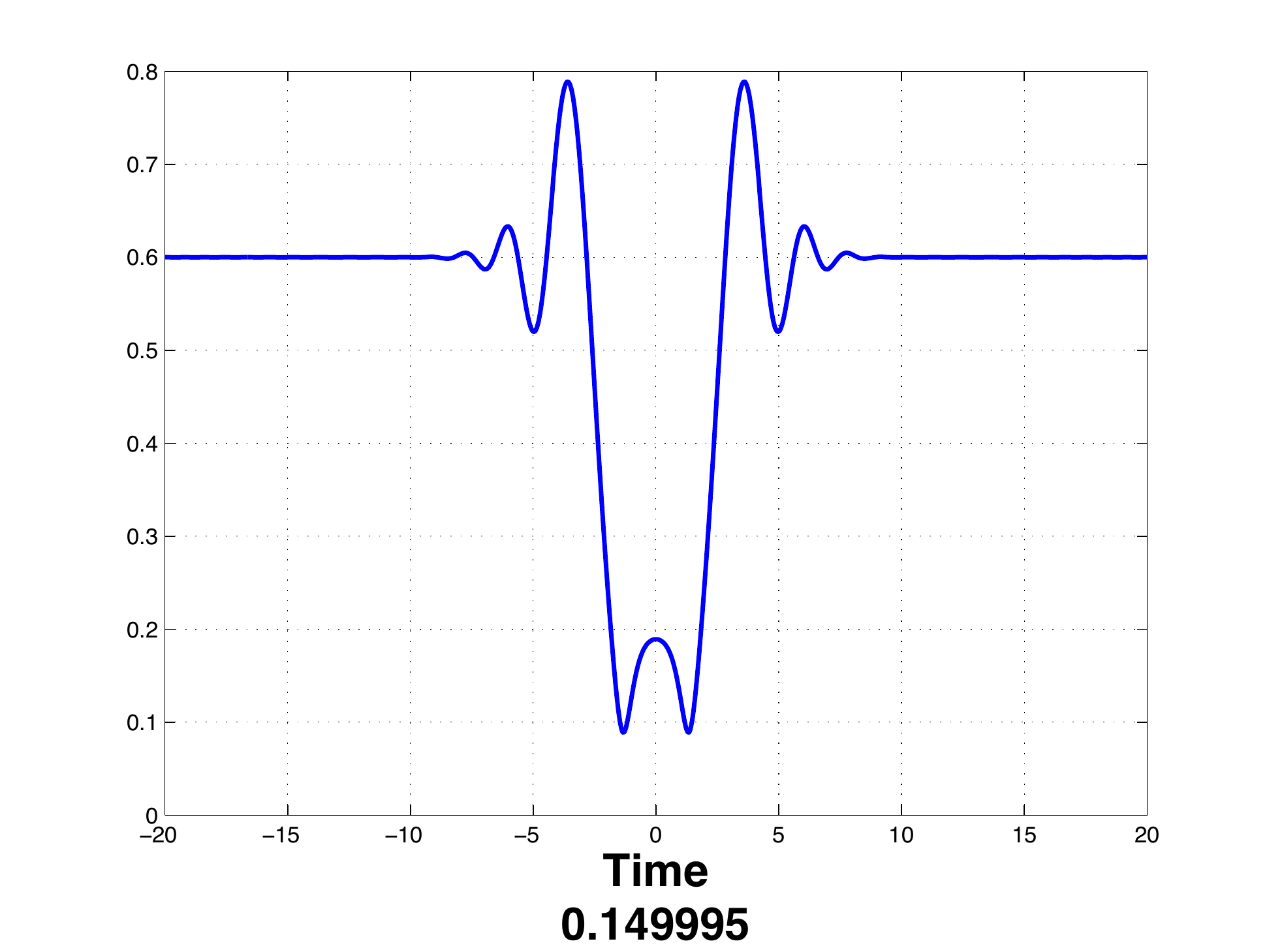}
   \includegraphics[height=4.2cm]{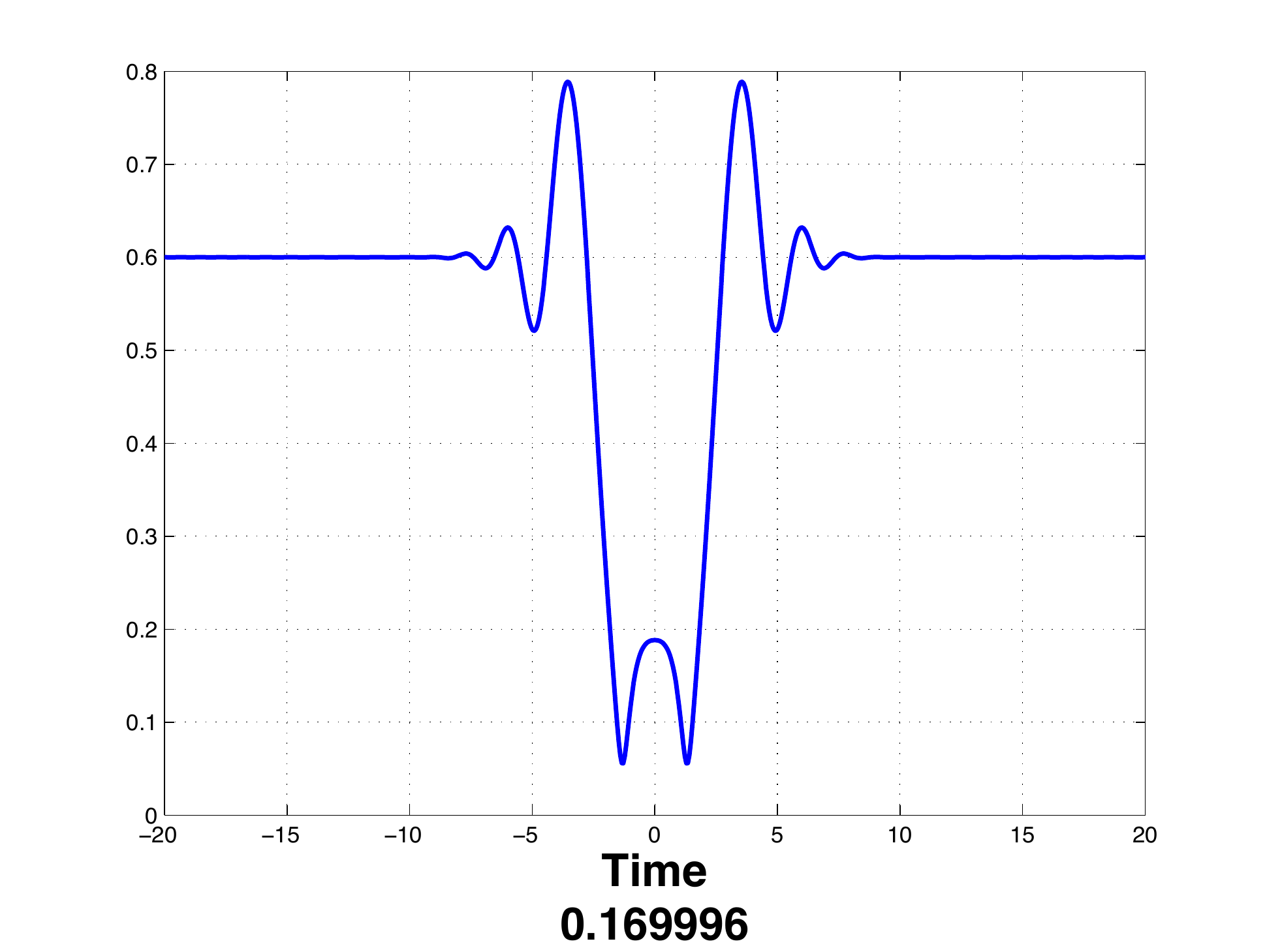} 
   \includegraphics[height=4.2cm]{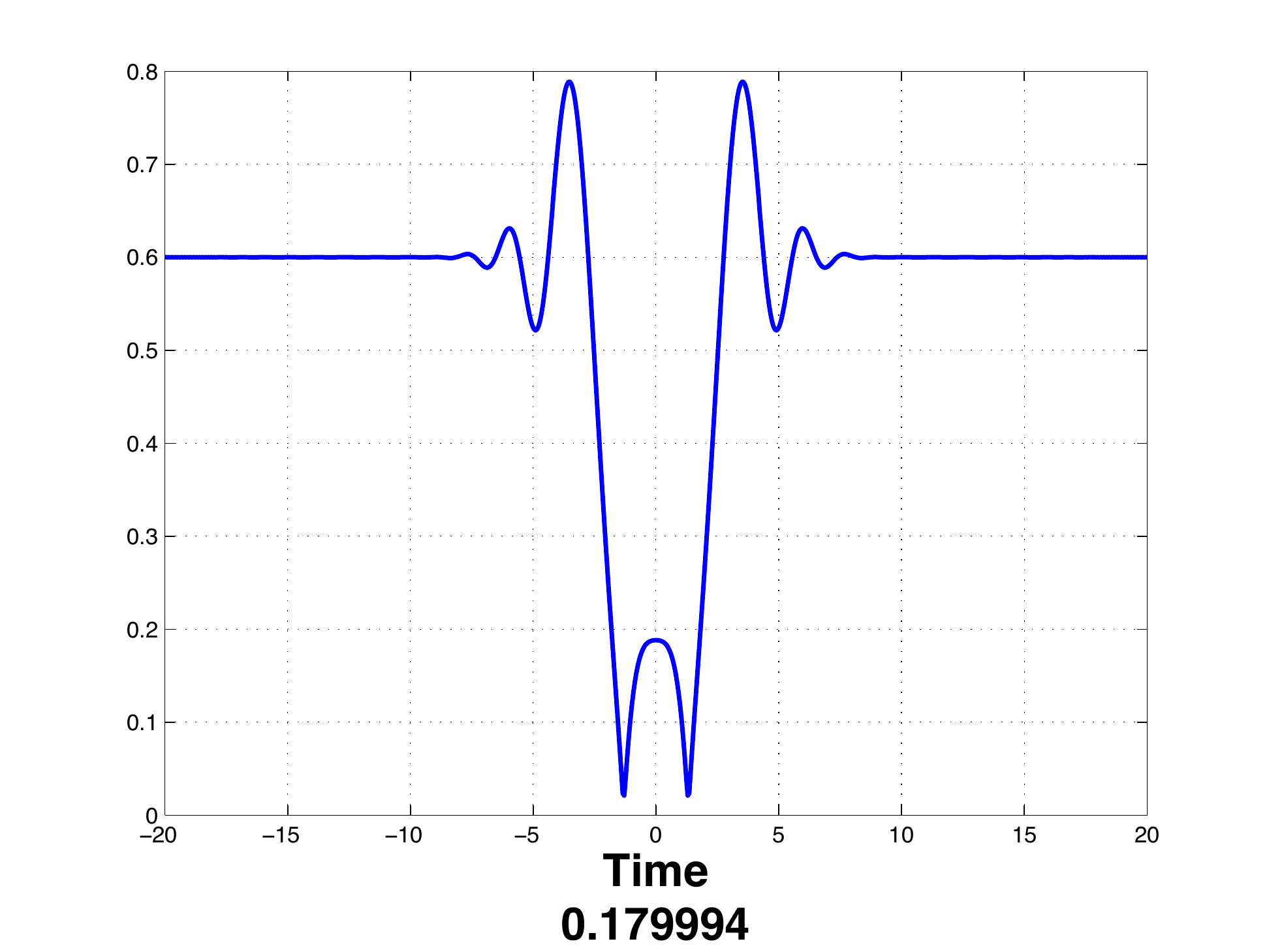},
    \includegraphics[height=4.2cm]{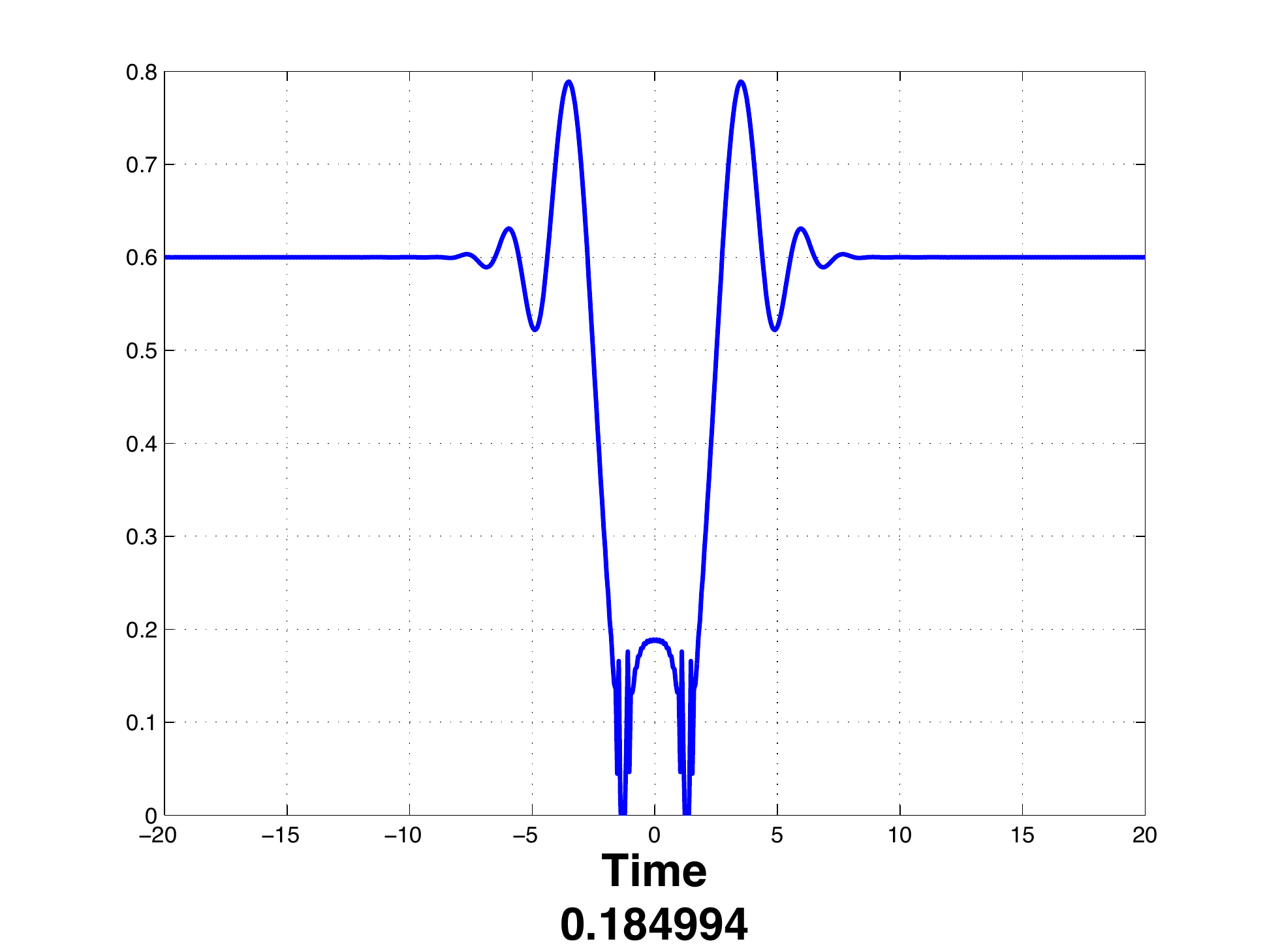}
  
  \caption{Real part of $\psi_1(t)$ for collision in the case  $\alpha_1=\alpha_2=1$, $\Psi_1(0)=0.6-G$, $\Psi_2(0)=-0.6+\overline{G}$, symmetry $\Psi_1(t)=-\overline{\Psi}_2(t)$.} \label{re2} \end{figure}

In the last figure  (see Figure \ref{re2}), we perform a numerical simulation of \eqref{eqpsi} but with initial data $\Psi_{1}(0) = 0.6 - G$. In this situation we observe a collision at a time close to $t = 0.18$. Note that in this situation, the collision arises at two symmetric points near $\sigma = 0$ We mention that this kind of collision is also observed in \cite{Ke93} where numerical simulations of the interaction of perturbed antiparallel vortex tubes are performed in the setting of three-dimensional incompressible Euler equations.  
We would like to conclude on this case by doing the following remarks. First, the collisions numerically observed are locally of a form very similar to the solutions constructed in the previous Section \S \ref{sect-selfsim}: at the collision time, the solution is locally made of two straight lines forming an angle on a collisional corner. Second, the numerical solution remains very close to the solution of the linear flow. Hence, by combining the two arguments (pertubative approach and description of a corner collision) these observations give a hope to rigorously describe the collision and prove its existence mathematically. 

\begin{thebibliography}{alpha}


\bibitem{BFG}
{\rm D. Bambusi, E. Faou and B. Gr\'ebert}
{\em Existence and stability of ground states for fully discrete approximations of the nonlinear Schr\"odinger equation}. Numer. Math. \textbf{123} (2013) 461--492

\bibitem[BM12]{BaMi} V. Banica and E. Miot, 
 {\it Global existence and collisions for symmetric configurations of nearly parallel vortex filaments}, 
 Ann. Inst. H. Poincar\'e Anal. Non Lin\'eaire \textbf{29} (2012), 813--832.

\bibitem[BM13]{BaMi2} V. Banica and E. Miot, 
{\it Evolution, interaction and collisions of vortex filaments},
Differential and Integral Equations \textbf{26} (2013), 355--388. 



\bibitem[BPSS13]{BPSS} J. L. Bona, G. Ponce, J.-C. Saut and C. Sparber,
{\it Dispersive blow up for nonlinear Schr\"odinger equations revisited}, ArXiv:1309.5023.



\bibitem[CG13]{CrGa} W. Craig and C. Garc\`ia-Azpeitia, 
\url{ http://www.math.uzh.ch/nhpde12/fileadmin/nhpde12/\\pdf/Craig\_Monday1030Ascona2012.pdf}.

\bibitem[C70]{Cr} S. C. Crow,
{\it Stability theory for a pair of trailing vortices}, 
AIAA J. \textbf{8} (1970), 2172--2179.

\bibitem[F12]{F11}
E. Faou,
{\em Geometric numerical integration and Schr\"odinger equations}, European Math. Soc., 2012. 

\bibitem[K93]{Ke93}
R. M. Kerr, 
{\it Evidence for a singularity of the three-dimensional, incompressible Euler equations}, 
Phys. Fluids A \textbf{5} (1993), 1725--1746.



\bibitem[KMD95]{KlMaDa}
R. Klein, A. J. Majda, and K. Damodaran,
 {\it Simplified equations for the interaction of nearly parallel vortex filaments},
Journal of Fluid Mechanics \textbf{288} (1995), 201--248.

\bibitem[KPV03]{KePoVe}
C. Kenig, G. Ponce, and L. Vega,
 {\it On the interaction of nearly parallel vortex filaments}, 
Comm. Math. Phys. \textbf{243} (2003), 471--483.

\bibitem[LM00]{LiMa} P.-L. Lions and A. J. Majda, 
{\it Equilibrium statistical theory for nearly parallel vortex filaments},
Comm. Pure Appl. Math. \textbf{53} (2000), 76--142.


\bibitem[MB02]{MaBe} A. J. Majda and A. L. Bertozzi, 
 "Vorticity and incompressible flow",
Cambridge Texts in Applied Mathematics, 2002.

\bibitem[MZ97]{MeZa}
F. Merle and H. Zaag, 
{\it Reconnection of vortex with the boundary and finite time quenching}, 
Nonlinearity \textbf{10} (1997), 1497--1550.

\bibitem[MZ97-2]{MeZa2}
F. Merle and H. Zaag, 
{\it Stability of the blow-up profile for equations of the type $u\sb t=\Delta u+\vert u\vert \sp {p-1}u$}, 
Duke Math. J. \textbf{86} (1997), 143--195.


\bibitem[Z88]{Zh1} V. E. Zakharov, {\it Wave collapse}, Sov. Phys. Usp. \textbf{31} (1988), no. 7, 672--674.

\bibitem[Z99]{Zh2} V. E. Zakharov, {\it Quasi-two-dimensional hydrodynamics and interaction of vortex tubes}, Lecture Notes in Physics \textbf{536} (1999), 369--385.



\end{thebibliography}
\end{document}